\input amstex
\documentstyle{amsppt}
\magnification=1200
\NoBlackBoxes
\TagsOnRight
\NoRunningHeads
\vsize=22 truecm
\hsize=16 truecm

\def \real{{\Bbb R}}
\def \complex{{\Bbb C}}

\def \sp {{ \text{sp}\, }}
\def \Id {{\text{Id}\, }}
\def \Fix {{\text{Fix}\, }}
\def \Det {{\text{Det}\, }}
\def \Ker {{\text{Ker}\, }}
\def  \Imm {{\text{Im}\, }}
\def \tr {{\text{tr}\, }}

\def \sgn{{\text{sgn} \, }}
\redefine\sp{\operatorname{sp}}

\def\AA{{\Cal A}}
\def\BB{{\Cal B}}

\def\DD{{\Cal D}}
\def\EE{{\Cal E}}
\def\FF{{\Cal F}}

\def\JJ{{\Cal J}}
\def\KK{{\Cal K}}
\def\LL{{\Cal L}}
\def\MM{{\Cal M}}
\def\NN{{\Cal N}}

\def\PP{{\Cal P}}
\def\QQ{{\Cal Q}}
\def\RR{{\Cal R}}
\def\SS{{\Cal S}}
\def\TT{{\Cal T}}

\def\today {\ifcase\month\or January \or February \or March \or
April \or May \or June
\or July \or August \or September \or October \or November \or December
\fi
\number\day~\number\year}

\document
\topmatter

\title\nofrills
Kneading determinants and spectra of transfer operators
in higher dimensions,\\
the isotropic  case
\endtitle

\author
Mathieu Baillif and Viviane Baladi
\endauthor

\date
June 2004 (Revised version)
\enddate

\address
V. Baladi: C.N.R.S.,
I.H.\'E.S.,  route de Chartres, 91440 Bures-sur-Yvette, FRANCE
\hfill\break
Current address: C.N.R.S., Institut Ma\-th\'e\-ma\-ti\-que de Jussieu,
F-75251 Paris, FRANCE
\endaddress
\email
baladi\@math.jussieu.fr
\endemail

\address
M. Baillif: Section de math\'ematiques, CH-1211 Gen\`eve 24, SUISSE
\endaddress
\email
baillif\@math.unige.ch
\endemail

\abstract
Transfer operators $\MM_k$ 
are associated to $C^{r}$
transversal local diffeomorphisms $\psi_\omega$
of $\real^n$, and $C^{r}$ compactly
supported functions $g_\omega$. A formal trace
$\tr^\# \MM$,  yields a
formal Ruelle-Lefschetz determinant $\Det^\#(\Id -z\MM)$.
We use the Milnor-Thurston-Kitaev equality recently proved by
Baillif [Bai]
to relate zeroes and poles  
of $\Det^\#(\Id -z\MM)$ with  spectra
of the transfer operators $\MM_k$,  under  additional assumptions.
As an application, we get a new proof of a result of Ruelle [Ru3] on the
spectral interpretation of zeroes and poles of the
dynamical zeta function 
$\exp \sum_{m\ge1} (z^m/m)
\sum_{f^m (x)=x} |\det Df(x)|^{-1}$
for smooth
expanding endomorphisms $f$.
\endabstract


\thanks
We thank
A.~ Avila, M. Benedicks, J.~ Buzzi, G.~ David, 
C. ~Liverani,   and D.~ Ruelle for useful conversations,
and S. ~ Gou\"ezel who found some errors in
previous versions of this article. We are grateful  to the referee for a careful reading
which helped us to improve the paper.
We benefitted from funding
by PRODYN and from the friendly environment of
l'Odyss\'ee dynamique 2001 at CIRM.
This work was started  
at the Universit\'e de Paris-Sud,  and continued at I.H.\'E.S. We express our thanks to both
institutions.
\endthanks

\endtopmatter

\head
1. Introduction
\endhead

A weighted (Ruelle) dynamical zeta function may be associated to a transformation
$f$ (on a compact space $M$, say) and a  ``weight'' function
$g:M\to \complex$ by setting
$$
\zeta_{f,g}(z) = \exp \sum_{m=1}^\infty {z^m \over m}
\sum_{x\in \Fix(f^m)} \prod_{j=0}^{m-1} g(f^j(x))\, ,\tag{1.1}
$$
assuming, e.g., that each set
$\Fix f^m = \{x \in M \mid f^m(x)=x \}$ is finite.
If $f$ is a one-sided subshift of finite type associated to a
Markov matrix $A$, and $g$ is constant, it is an easy algebraic exercise to see that
$\zeta_{f,g}(z)=1/\det(\Id -zA)$. 
If $g$ is not constant, but has exponentially decaying variations
on small cylinders, the transfer operator
$$
\LL \varphi(x)=\sum_{f(y)=x} g(y) \varphi(y) \tag{1.2}
$$
has essential spectral radius bounded above by $\theta\exp(P\log(|g|))$
(with $P$ the topological pressure) when acting on the Banach space
$\text{Lip}$ of Lipschitz functions on the one-sided shift space
([Po]).
(Throughout this introduction, we refer to
the survey [Ba1] for more information and references.)
Haydn then proved (following crucial
steps by Ruelle and Pollicott) that $\zeta_{f,g}(z)$ is meromorphic
in the disc of radius  $\theta^{-1}\exp(-P\log(|g|))$, where its
poles are the inverse eigenvalues of $\LL$ on  $\text{Lip}$. 

Using Markov methods
(symbolic dynamics via
Markov partitions, Markov covers, tower-like Markov extensions, etc.), similar results were obtained
for weighted dynamical zeta functions associated to smooth
data on manifolds: Axiom A
diffeomorphisms (Ruelle [Ru2], Pollicott [Po], Haydn [Ha]), one-dimensional piecewise monotone maps
(Hofbauer-Keller, Baladi-Keller [BaKe], Keller-Nowicki [KN]), and higher-dimensional
piecewise smooth and injective maps (Buzzi-Keller [BuKe]).
There, one assumes not only (piecewise) smoothness of the dynamics $f$ and
the weight $g$ but also some ``hyperbolicity'' (or expansion) in $f$.

Note also that, for smooth expanding
maps on manifolds, Ruelle [Ru1, Ru3] and Fried [Fr0, Fr1] studied dynamical Fredholm-type
determinants. These determinants
will only
be mentioned in passing in the present article (see \thetag{4.2} below). 
We just point out here that, in the analytic case, 
the transfer operator is nuclear on
a suitable Banach space, and this dynamical determinant can further
be interpreted as an honest Grothendieck-Fredholm determinant:
see  Mayer [Ma], Ruelle [Ru1], Keller [Ke], 
and more recently Rugh [Rug], Fried [Fr1]. The underlying topological
mechanism in these works is still the Markov partition.

\smallskip

In the early nineties, a  new approach was launched [BaRu and references therein],
initially for piecewise monotone interval maps.
The original motivation was to understand the links between, on the one hand the
Milnor-Thurston  [MT] identity relating the kneading matrix 
and an unweighted dynamical zeta function, and 
on the other hand the spectral interpretation of the zeroes and
poles of (weighted) dynamical zeta functions.
This kneading approach is now quite well understood in
one real dimension
([Ru4, BaRu, BKRS, Ba2, Go]),
We recall its key ingredients next, very briefly.

The initial data
is a finite set of local homeomorphisms $\psi_\omega: U_\omega \to \psi_\omega(U_\omega)$,
where each $U_\omega$ is
an open interval of
$\Bbb R$, and
of associated weight functions $g_\omega$
which are continuous,
of bounded variation,  and have support inside 
$U_\omega$. In particular, the $\psi_\omega$ can be the inverse branches of
a single piecewise monotone interval map $f$, and $g_\omega$ can be
$g\circ \psi_\omega$ for a single $g$. For a more specific
example, if $f$ is unimodal
on $[-1,1]$ with $\max f=f(0)=b<1$ and $\min f=f(+1)=f(-1)=-1$
then we may take $U_-=U_+=(b+\epsilon,-1-\epsilon)$, and
$\psi_\pm$  (somewhat arbitrary
extensions of) the two inverse branches of $f$, with images in
$(-1-\delta, \delta)$ and $(-\delta, 1+\delta)$. We may also
take $g$ a continuous function
of bounded variation on $(-1-\delta, 1+\delta)$ supported in $[-1,+1]$.
Beware however that
this setup does {\it not} allow to take for $g$ the characteristic
function of $[-1,1]$.
Gou\"ezel [Go]  recently obtained a
significant weakening of the continuity assumption 
(while avoiding boundary terms in integration by parts),
which allows e.g. the choice $g=1_{[-1,+1]}$ 
for our unimodal example {\it if} $0$ is not preperiodic, i.e., if there is no finite Markov
partition. 

We emphasize that no contraction assumption is required on the $\psi_\omega$:
their graph can even coincide with the diagonal on a segment.
The transfer operator is now 
$$
\MM \varphi=\sum_\omega g_\omega \cdot (\varphi \circ \psi_\omega) \, .\tag{1.3}
$$
Ruelle [Ru4] obtained an estimate, noted $\widehat R$, for the
essential spectral radius of $\MM$ acting on
the Banach space $BV$ of functions of bounded variation.
The main result of [BaRu] links the eigenvalues of $\MM$ 
(outside of the disc of radius $\widehat R$),  acting
on $BV$, with the zeroes of the ``sharp determinant''
$$
\Det^\# (\Id-z\MM)=\exp- \sum_{m=1}^\infty
{z^m \over m} \tr^\# \MM^m \, , \tag{1.4}
$$
where (with the understanding that $y/|y|=0$ if $y=0$)
$$
\tr^\# \MM =   \sum_\omega \int {1\over 2}  {\psi_\omega(x)-x\over
|\psi_\omega(x)-x|} \, dg_\omega(x) \, . \tag{1.5}
$$
If the $\psi_\omega$ are strict contractions which form the set of
inverse branches of a piecewise monotone interval map $f$, and
$g_\omega=g \circ \psi_\omega$ then integration
by parts together with the key property that
$$
d {x \over 2 |x|}      =\delta\, , \hbox{ the Dirac delta at the origin of }   \real\, ,
$$
show that
$\Det^\# (\Id -z\MM)=1/\zeta_{f,g}(z)$. If one only assumes that
the graph of each admissible composition
$\psi^m_{\vec \omega}$ of $m$ successive $\psi_\omega$s (with $m\ge 1$)
intersects the diagonal transversally, then
$$
\Det^\#(\Id -z\MM) =\exp -\sum_{m=1}^\infty
{z^m \over m} \sum_{\hbox{admissible } \psi^m_{\vec \omega}}
\sum_{x \in \Fix \psi^m_{\vec \omega}} L(x,\psi^m_{\vec \omega})
\prod_{j=0}^{m-1} g_{\omega_j}(\psi^j_{\vec \omega}(x)) \, , \tag{1.6}
$$
where $L(x,\psi)\in \{-1,1\}$ is the Lefschetz number of a
transversal fixed point $x=\psi(x)$ (if $\psi$ is $C^1$
this is just $\sgn(1-\psi'(x))$). This is why we  call
the sharp determinant
$\Det^\#(\Id -z\MM)$ a Ruelle-Lefschetz (dynamical) determinant.
It is not difficult to see,
for non Markov unimodal $f$ and $g=1_{[-1,1]}$
as above, that the expression \thetag{1.6} with Lefschetz numbers,
coming from our additional
transversality assumption, gives that
$\Det^\#(\Id -z\MM)$ is just $1/\zeta^{-}(z)$, where
the ``negative zeta function''
$$
\zeta^{-}(z)=\exp +\sum_{m=1}^\infty
{z^m\over m} \bigl (2 \# \Fix^{-} (f^m) -1 \bigr )
$$
is defined by counting the sets
$$ 
\Fix^{-} (f^m)=\{x \in (-1,1]\mid f^m (x)=x \, ,
f \hbox{ strictly decreasing in a nbhd of } x\}
$$ 
of ``negative fixed points.'' This zeta function was studied by Milnor and Thurston [MT],
who proved the remarkable identity
$$
(\zeta^{-}(z))^{-1} = \det (1+\widehat D(z))
$$
where $\widehat D(z)$ is a  $1\times 1$ ``matrix'' which is a power series in $z$
with coefficients in $\{-1,0,+1\}$ given by the signed itinerary of $f(0)$
(i.e., the kneading data).
Milnor and Thurston's equation also applies to unimodal maps,
but $\widehat D(z)$ is here a finite $(N-1)\times( N-1)$ matrix (with $N$
the number of laps of $f$) constructed from the signed kneading series
of the $N-1$ turning points.

Let us return now to our general setup $\psi_\omega$, $g_\omega$.
The crucial step in the proof of the spectral interpretation
of the zeroes of this Ruelle-Lefschetz determinant consists in establishing
the following continuous version of the Milnor-Thurston
identity:
$$
\Det^\#(\Id -z\MM) = \Det^ *(\Id+\widehat \DD(z)) \, , \tag{1.7}
$$
where the ``kneading operator''
$\widehat \DD(z)$     replaces (formally)
the finite kneading matrix of Milnor and Thurston.
(In fact, approaching the $g_\omega$ by locally constant functions,
$\Det^*(\Id+\widehat \DD(z))$ can be approached by the determinants of
a sequence of finite weighted Milnor Thurston
matrices, see the appendix of [BaRu].)

If $|z|< 1/R$  then $\widehat \DD(z)$  is
a Hilbert-Schmidt operator on a suitable $L^2$ space
(its kernel is in fact bounded and compactly supported),
thus allowing the use of regularised determinants
of order two. In fact, $\Det^*(\Id+\widehat \DD(z))$ is the product of the regularised
determinant with the exponential of the average of the kernel of $\widehat \DD(z)$
along the diagonal, which is well-defined. Another kneading
operator, $\DD(z)$, is essential. If $1/z$ is not in the spectrum of $\MM$
(on $BV$) then $\DD(z)$  is also Hilbert-Schmidt, and we have
$\Det^ *(\Id+\widehat \DD(z))=\Det^ *(\Id+ \DD(z))^{-1}$. 

The kneading operator approach
to study Ruelle-Lefschetz dynamical determinants  in one
complex dimension is presented in [BKRS]. There,  a more conceptual 
definition of the $\DD(z)$ was suggested, which
was later implemented in one real dimension [Ba2]:
$$
\DD(z)=\NN (\Id-z\MM)^{-1} \SS \, , \tag{1.8}
$$ 
where $\NN$ is an
auxiliary transfer operator and $\SS$ is the convolution 
($\mu$ being an auxiliary nonnegative finite measure)
$$
\SS \varphi(x) = \int
{1 \over 2} {x-y\over |x-y|} \varphi(y) d\mu \, . \tag{1.9}
$$
It  is indeed clear from \thetag{1.8} 
that the kneading operator is a regularised (through the convolution $\SS$)
object which describes the inverse spectrum
of the transfer operator: the resolvent $ (\Id-z\MM)^{-1}$ in \thetag{1.8}
means that poles can only appear if $1/z$ is an eigenvalue,
and since $\Det^ *(\Id+\widehat \DD(z))=\Det^ *(\Id+ \DD(z))^{-1}$,
this can be translated into a statement for zeroes of
$\Det^ *(\Id+\widehat \DD(z))$.
The Milnor-Thurston
identity \thetag{1.7}) then implies that any zero
of  $\Det^\#(\Id -z\MM)$ (recall $|z|<1/R$) is an inverse eigenvalue
of $\MM$.

In our opinion, the kneading approach we just described
is not only interesting because
it extends previous  results (see e.g. [BaKe])
on piecewise monotone interval maps, 
but also because of its conceptual
simplicity. We believe it
sheds a new light on weighted dynamical zeta functions associated
to smooth (non analytic) maps: The mechanism relating the
zeroes of the zeta function with eigenvalues is 
literally visible in \thetag{1.7--1.8}, as we just explained. As a side-effect, the sometimes
cumbersome Markov partition tool is bypassed. 
(The unimodal example discussed above shows that in fact this
one-dimensional approach works best for piecewise monotone
maps with no finite Markov partition.)

\smallskip

Implementing this strategy in higher dimensions
is a natural goal: Assume that the $U_\omega$
are (finitely many, say) open subsets of $\real^n$ and that the $\psi_\omega:U_\omega
\to \psi_\omega(U_\omega)$
are local $C^r$ homeomorphisms or diffeomorphisms, while the
$g_\omega:U_\omega\to \complex$ are compactly
supported $C^r$ functions, for $r\ge 1$. 
We assume that $\cup_\omega U_\omega$ is contained
in a compact set. We consider
the $n$-tuple $\MM$ of operators $\MM_k$ acting on $k$-forms by
$$
\MM_k \varphi=\sum_\omega g_\omega \cdot \psi_\omega^*  \varphi \, .\tag{1.10}
$$
As in dimension one, constant (or discontinuous) weights  $g_\omega$ are
not allowed. Since we may reduce to the model above from a smooth dynamical
system on a manifold with the help of charts and partitions of unity, this
causes in fact less problems than in  dimension one.
For example the perhaps most natural  transfer operator acting on functions,
$$
\LL_f \varphi (x)=\sum_{f(y)} \varphi(y)/|\det D(f(y))|\, ,
$$
where $f$ is the expanding
linear map on the two
torus defined by $(y_1,y_2)\mapsto (2 y_1, 2y_2)$ (modulo $1$),
and $\det Df \equiv 4$, can easily be modelled by
\thetag{1.10} with the help of local inverse branches and partitions of unity. We refer
to Section ~4 for details. (If $f$ is an Anosov toral map, a partition of unity can
also be used to obtain a representation
\thetag{1.10}, without the need of Markov partitions.)

In 1995, A.~Kitaev wrote a two-page sketch
of the (statement and)
proof of a higher-dimensional Milnor-Thurston formula. He made an additional transversality
assumption, which guarantees that the set of fixed points
of each fixed period  $m$ is finite, allowing
him to define a Ruelle-Lefschetz determinant
$
\Det^\#(\Id -z \MM)$
through formula \thetag{1.6}.
(Note that the Lefschetz numbers $ L(x,\psi^m_{\vec \omega})
=\sgn( \det (\Id -D_x \psi^m_{\vec \omega}))$ in \thetag{1.6} 
are nonzero by the transversality assumption, and that if we assume that
the $\psi_\omega$ are contractions,  they are all $+1$.)
Baillif
turned  Kitaev's unpublished note into a mathematical
proof  [Bai, see also Ba2] and proved the following ``Milnor--Thurston" formula:
$$
\Det^\#(\Id -z \MM)
=\prod_{k=0}^{n-1} \Det^\flat(\Id+\DD_k(z))^{(-1)^{k+1}}    \, . \tag{1.11}
 $$
We refer the reader to Section ~2 for the definition of the kneading
operators $\DD_k(z)$ and their flat traces and determinants. We just mention here that they are
kernel operators acting on
$k+1$ forms,   constructed  with the resolvent
$(\Id-z\MM_k)^{-1}$ (where the transfer operator $\MM_k$
acts on $k$ forms via the pull-back $\psi_\omega^*$), together with a convolution 
operator $\SS_k$, mapping
$k+1$-forms to $k$-forms and
which satisfies the homotopy equation $d\SS+\SS d=1$. The
kernel $\sigma_k(x,y)$ of $\SS_k$  has singularities
of the form $(x-y)/\|x-y\|^n$.
The transversality assumption allows Baillif
to interpret the determinant obtained by integrating
the kernels along the diagonal as a flat determinant
in the sense of Atiyah and Bott [AB1, AB2], whence the notation $\Det^\flat$
in the right-hand-side of \thetag{1.11}.

Although Baillif did not give a spectral interpretation of zeroes
or poles of the sharp determinant \thetag{1.11},
he did notice that for $|z|$ very small, suitably high iterates of the $\DD_k(z)$
are trace-class on $L^2(\real^n)$, showing that the corresponding
regularised determinant has a nonzero
radius of convergence under rather weak assumptions.

\smallskip
In the present work, we carry out the spectral interpretation
of the sharp determinant \thetag{1.6} in arbitrary dimension, under the transversality
assumptions of Baillif and Kitaev (see Section~3).
We concentrate here on the isotropic setting where the
forward dynamics is (uniformly) expanding in all space directions.
We decided however to formulate our main result in a way suitable
for generalisations to anisotropic (in particular, Anosov) cases.
This result says   (Theorems~13--14 in Section~3) that
{\it if} (Axiom~ 1 in Section~3)
there are suitable  Banach spaces $\BB_k$ of coefficients
for $k$-forms on which (Axiom ~2 in Section~3)
the essential spectral radii of
all transfer operators $\MM_k$ are smaller than  some $\widetilde R<1$,
while the flat traces of the kernels $\DD_k(z)$ are meromorphic in
the disc of radius $1/R$ (Axiom~3), then
{\it taking $R=\widetilde R^{1/n}$:}

\proclaim{Main Theorem} Under Axioms 1, 2, 3 from Section~3,
$\Det^\#(\Id-z\MM)$ is meromorphic in $\{ |z| < 1/  R\}$.
The order of $z$ as a zero/pole of $\Det^\#(\Id-z\MM)$ is

\noindent $\bullet$
the sum over $0 \le 2k \le n$ of the algebraic multiplicity
of $1/z$ as an eigenvalue of   $\MM_{2k}$ acting on $2k$ forms with
coefficients in $\BB_{2k}$ 

\centerline{\bf minus} 

\noindent $\bullet$ the sum over
$1 \le 2k+1 \le n$ of the algebraic multiplicity
of $1/z$ as an eigenvalue of  $\MM_{2k+1}$ acting on $2k+1$-forms with
coefficients in $\BB_{2k+1}$. 
\endproclaim

When the dynamics
$\psi_\omega$ corresponds to the inverse branches of a differentiable uniformly
expanding map on a compact manifold and the $g_\omega$ are differentiable,
Axioms ~1 and ~2--3 are satisfied by taking
the $\BB_k$ to be Sobolev spaces.  We then  recover 
a (case of a) result of Ruelle [Ru3] (Section~ 4).

\smallskip

{\bf Choosing  Banach spaces}

\smallskip
In dimension one,  the kernel  of the convolution operator $\SS$ is a function of
bounded variation. This is very handy,
because  the essential spectral radius of the transfer operator
on $BV$ can be estimated, and in many cases shown to be
strictly smaller than the spectral radius (``quasicompactness'').
To adapt the strategy of [BaRu] to our higher-dimensional setting, it is
desirable to find Banach spaces which, on the one hand,
contain the kernels $\sigma_k(\cdot, y)$ (or some ``primitive'' of the
kernels, see below) and,
on the other, allow us to obtain quasicompactness of
the transfer operators.

The kernels $\sigma_k(\cdot,y)$  belong to a fractional
Sobolev space $W^{\alpha,p}$ for $0< \alpha < 1$ and small
enough $p > 1$. Instead of working with this
limited ``H\"older''-type fractional regularity,
we adapt 
a technique used by Ruelle [Ru4] in dimension one:
he constructed modified kneading operators $\DD^{(r)}(z)$
with smoother
kernels obtained by taking a ``primitive'' of the $\sigma_k$, which have the same traces as the original
operators. This step is carried out in the abstract Lemma ~6.
In the concrete setting of Section 4, it allows us to let the transfer operators
act on forms  whose coefficients are (locally supported and)
in a Sobolev space $L^{t}_m(\real^n)=W^{m,t}(\real^n)$, i.e.,  have their
derivatives in  $L^{t}(\real^n)$, up to some finite order $m\ge 0$.
To achieve this goal, we need an analogue of the one-dimensional
concept of ``primitive.'' An obvious tool [St] is the inverse square
root of the Laplacian (Riesz potential), or its better behaved
version at infinity, the Bessel potential.  
This step is emblematic of the most important new feature
of this higher-dimensional version of [BaRu]: the {\it intervention
of harmonic analysis.}  Indeed, very classical
Calder\'on-Zygmund results are behind the fact (Lemma~3) that
the operators $d\SS_k$ are bounded on $L^p(\real^n)$.
The second new ingredient is the (re)regularisation \thetag{3.9}
of the kneading operators $\DD_k(z)$ (in fact; of their smoothened
version $\DD^{(r)}_k(z)$): In particular,
the resolvent  $(1-z\MM_k)^{-1}$ is replaced by  (see \thetag{3.6})
$(1-\MM_{k,L}(z))^{-1}$, where
$
\MM_{k,L}(z) =\Id - (\Id - z\MM_k)
\exp\biggl (\sum_{j=1}^{L-1}{z^j\MM_k^j\over j} \biggr ) 
$,
for large enough $L$. (This allows us to ``kill'' constants
and essentially consider the kneading operators as compact,
see Corollary of Lemma~5.)
Thirdly,  one of the difficulties of the higher-dimensional kneading approach
is that $d$ and $\MM_k$ do not commute.
The spectra  of operators $d\SS_{k-1} \MM_k$ 
(note that $d\SS$ is a projection to $\Ker d$) arise naturally from
the definitions of the kneading operators (see Lemma~11).  In order to show 
that, outside of the disc of radius $R$, the  eigenvalues $\lambda$ of $d\SS_{k-1}\MM_k$ 
which are not
eigenvalues of $\MM_k$  do not affect the zeroes and
poles of the sharp determinant (in the corresponding disc), we show that
they are not intrinsic (see Lemma~12).  An indirect argument
(Lemma ~11, Theorem~ 13) then proves that the contributions
of such eigenvalues must cancel
in the alternated product of kneading determinants.
These two other new tricks are simple-minded,
but they seem  not to have been used before in dynamical contexts.
We hope they will be helpful elsewhere, the last
one  may be relevant in dynamical cohomological settings.

\smallskip

{\bf Further extensions -- beyond the smooth isotropic case}

\smallskip

A natural extension is to apply
kneading operators to uniformly hyperbolic dynamics, in order to  improve on the
results  of Pollicott, Ruelle, and Haydn mentioned above 
and to   contribute to a spectral interpretation of the
work of Kitaev [Ki], see  Fried [Fr2] and Blank-Keller-Liverani
[BKL]. In a forthcoming work by the second named author, we shall see that replacing the isotropic
Sobolev spaces of Section~4 
by   Banach spaces of anisotropic 
distributions, 
with leafwise smoothness, is the key to this extension. 
We also hope that the methods introduced in the present work will
eventually be used to tackle the nonuniformly (expanding or) hyperbolic case
(such as H\'enon-like maps).
A less ambitious but  interesting problem is to adapt these methods
to  higher-dimensional (non-Markov)
piecewise expanding maps, in order to give a different
proof (and perhaps improve on) the results of  Buzzi and Keller [BuKe].
Ultimately, one would also like to get rid of the transversality assumption.

\smallskip

{\bf Contents of the paper}

\smallskip
In Section 2 we introduce notations, recall the Kitaev-Milnor-Thurston
formula and prove some facts about the operators $\SS_k$.
Section 3 contains a statement of our
assumptions (Axioms 1 and 2--3) and our main result (Theorem 13).
(The reader will notice that the axioms below are more general 
than the ones stated in [Ba3], but we have to take
$R=\widetilde R^{1/n}$. We have also corrected 
Axioms 2 and~3 there.) Assuming only Axioms~1 and 2, we get
a weaker result (Theorem 14).
In Section 4 we apply Theorem 13 to the smooth expanding setting, 
recovering a result of Ruelle ~[Ru3].


\smallskip

\head
2. Preliminaries -- the Milnor and Thurston formula
\endhead

{\bf Notations and basic data - $\psi_\omega$ and $g_\omega$}
\smallskip

Let us fix once and for all an integer order of differentiability
$r \ge 1$ and a dimension $n\ge 2$. 
The case $n=1$ is dealt with in [BaRu, Ru4] and [Go] (see also [Ba1, Ba2]).

The data we consider is a  system (indexed
by a finite set $\Omega$) of
\roster
\item
local $C^r$ diffeomorphisms
$\psi_\omega : U_\omega \to \psi_\omega(U_\omega)$, where
$r\ge 1$ is an integer, and each
$U_\omega$ is a nonempty bounded open subset of $\real^n$;
\item
local $C^r$ ``weight functions'' $g_\omega : \real^n
\to \complex$, such that the support of $g_\omega$ is
contained in $U_\omega$.
\endroster

We need more hypotheses and notation.  Let
$K = \{ \| x \| \le T \}$ in $\real^n$ be so that
$\cup_\omega U_\omega \subset K$ and
$\cup_\omega \psi_\omega(U_\omega) \subset K$.
We make the following {\it transversality}
assumption: for each $x\in \real^n$ and $m\ge 1$ such that there
are $\omega_1, \ldots, \omega_m$ for which the
composition
$$
\psi^m_{\vec \omega}=
\psi_{\omega_m} \circ \cdots \circ \psi_{\omega_1} \tag{2.1}
$$
is well defined in a neighbourhood of $x$ and fixes
$x$, the derivative 
$
D_x (\psi_{\omega_m} \circ \cdots \circ \psi_{\omega_1})
$
does not have $1$ as an eigenvalue.
We  often consider products of weights along orbits
$$
g^{(\ell)}_{\vec \omega}(x):=
g_{\omega_\ell}(\psi_{\omega_{\ell-1}} \circ \cdots \circ \psi_{\omega_1} (x))
 \cdots g_{\omega_2} (\psi_{\omega_1}(x)) \cdot g_{\omega_1}(x)
 \, . \tag{2.2}
$$
Since the $m$th factor
$g_{\omega_m}(\psi_{\omega_{m-1}}  \cdots\psi_{\omega_1} (x))$
in the above product vanishes whenever the
composition \thetag{2.1} is not well defined (assuming the previous
compositions, for $j < m$ make sense), we shall not
worry much about domains of definition.
Note also that the transversality hypothesis,
combined with our compactness assumption, implies that for
each $\ell$, the
set $\Fix \psi^\ell_{\vec \omega}$ contains only finitely many points $x$
for which \thetag{2.2} isn't zero.

\smallskip

{\bf The transfer operators $\MM_k$ -- Forms $\AA_k$ -- Sharp determinant}

\smallskip

We write $\AA_{k,C^m}$, $\AA_{k,C^m(K)}$ for the vector spaces
of $k$-forms on $\real^n$ with $C^m$ coefficients ($m \ge 0$ integer),
respectively $C^m$ coefficients supported in the compact
set $K$. We shall also work with
$k$-forms with coefficients in $L^q=L^q(\hbox{Lebesgue}\,, \real^n)$
or in $L^q(U)=L^q(\hbox{Lebesgue}\, ,U)$, with
$U$ a bounded domain in $\real$
and $1 \le q < \infty$, writing $\AA_{k, L^q}$,
$\AA_{k, L^q(U)}$ for the
corresponding spaces.
More generally, if $\BB$ is a Banach space of functions
or distributions on $\real^n$, then $\AA_{k,\BB}$ denotes the space of $k$-forms with
coefficients in $\BB$. Sometimes we also consider
$\AA_{k,\BB}^C$    the space of compactly supported $k$-forms
with coefficients in $\BB$.
If $\BB$ is a Banach space, we use the
Banach norm 
$
\| \phi\|_{\BB} =\| \phi\|_{\AA_{k,\BB}}=\max_ {j\in I(k)} \|\phi_j\|_\BB 
$
for   $\phi  \in \AA_{k,\BB}$,
where $I(k)$ denotes the set of ordered $k$-tuples in
$\{1, \ldots, n\}$ and
$$
\phi(x)=\sum_{j\in I(k)} \phi_j(x) \, dx_{j_1} \wedge \cdots \wedge dx_{j_k} \, ,
\hbox{ with }\phi_j \in \BB   \, .
$$

The map $\psi_\omega$
induces a transformation $\psi_\omega^*$ on $k$-forms
by the usual pullback operation (which involves differentiating
the $\psi_\omega$ if $k\ne 0$).
Our object of interest is the transfer
operator associated to the system $\{ \psi_\omega, g_\omega\}$
by setting, for $k=0, \ldots, n$,
$$
\MM_k : \AA_{k,C^{r-1}} \to \AA_{k,C^{r-1}(K)} \, , \qquad
\MM_k  \phi(x) = \sum_\omega g_\omega(x) \, (\psi_\omega^* \phi) (x) \, .
\tag{2.3}
$$
(In fact, $\MM_0$ maps $\AA_{0, C^r}$ to $\AA_{0,C^r(K)}$.)
We write $\MM$ for  the $(n+1)$-tuple $(\MM_0, \ldots, \MM_{n})$.
We want to relate the spectra of the
$\MM_k$ on suitable Banach spaces  $\AA_{k, \BB_k}$ to the zeroes and poles
of the sharp determinant [Bai] associated to
the data $\{\psi_\omega, g_\omega\}$ by the (a priori) formal
power series
$$
\Det^\# (\Id-z\MM)=
\exp -\sum_{m=1}^\infty {z^m\over m} \sum_{\vec \omega\in \Omega^m}
 \sum_{x\in
\Fix \psi^m_{\vec \omega}}
g^{(m)}_{\vec \omega}(x)
\, L(x, \psi^m_{\vec \omega}) \, , \tag{2.4}
$$
where $L(x, \psi)\in \{ 1, -1\}$ (the value $0$ is excluded
by our transversality assumption) is the Lefschetz number of the diffeomorphism
$\psi$ at $x$, i.e.
$
L(x,\psi) = \sgn \det (1-D_x \psi) 
$.

In other words $\Det^\# (\Id-z\MM)=\exp -\sum_{m=1}^\infty {z^m\over m}\tr^\# \MM^m$,
the (formal) sharp trace being defined by
$
\tr^\# \MM = \sum_{\omega\in \Omega} \sum_{x\in\Fix \psi_\omega} 
g_{\omega}(x) \, L(x, \psi_{\omega}) 
$.

\smallskip

{\bf Kneading operators $\DD_k(z)$ -- Homotopy operators $\SS_k$}
\smallskip

We  shall introduce  next the
kneading operators
$\DD_k(z)$, for $k=0, \ldots, n-1$. Each  $\DD_k (z)$ is
a power series in $z$ with coefficients operators $\AA_{k+1,C^{r-1}(K)} \to
\AA_{k+1,C^{r-1}(K)}$ (other Banach
spaces of coefficients will be specified later).
For some values of $z$, these operators will have
a kernel and, for a suitable iterate $\DD^{m(n)}_k(z)$,
this kernel will be in $L^2$, so that
Hilbert-Schmidt theory will furnish  a regularised determinant
of order $2m(n)$.
Note that we
use the sign conventions of [Bai] when integrating forms depending on two
variables, for instance when considering kernel operators. For example
if $\phi(x,y)$ is a $n$-form in $y$, this
definition  implies 
$d^x\int_y\phi(x,y) = (-1)^n\int d^x\phi(x,y)$, where
$d^x$ denotes derivative with respect to $x$ (see  [Bai (5)--(7)])

The operators $\DD_k(z)$
were first defined in [Bai], following an idea of Kitaev.
They involve
the auxiliary transfer operators
$\NN_k : \AA_{k,C^{r-1}} \to \AA_{k+1,C^{r-1}(K)}$:
$$
\NN_k = d_k\MM_k - \MM_{k+1} d_{k}\, , \tag{2.5}
$$
where $d=d_k : \AA_{k,C^{r-1}} \to \AA_{k+1,C^{r-2}}$
(or $\AA^C_{k,C^{r-1}} \to \AA^C_{k+1, C^{r-2}}$)
denotes the exterior derivative. The fact that $\NN_k$
is well-defined for $r=1$ and does not lower the regularity comes,
e.g., from the equivalent expression
$
\NN_k \phi(x)=\sum_\omega dg_\omega(x) \wedge \psi_\omega^* \phi(x) 
$.

The good properties of the
kernels of the kneading operators
$\DD_k(z)$ are due to the 
convolution operators $\SS_k : \AA^C_{k+1,C^{r-1}} \to \AA_{k, C^{r-1}}$.
We shall see in the proof of Lemma~4 that $\SS_k$ can be written as
$d^*_{k+1} \Delta^{-1}_{k+1}$, we recall here the definition suggested
by Kitaev and used in [Bai]: Introduce an $n-1$ form solving $d\sigma = \delta$, the Dirac
current at  $0\in \real^n$,
$$
\sigma (x)  =
{\Gamma(n/2) \over 2 (\pi^{n/2})}\sum_{ i=1}^n \sigma^{\hat \imath} (x)
\, dx_1 \wedge \cdots
\wedge  d\hat x_i \wedge \cdots \wedge dx_n \, ,\,
\sigma^{\hat \imath} (x)=(-1)^{i+1} {x_i\over |x|^n}\, ,  \tag{2.6}
$$
where $\Gamma$ is the Euler function and the hat means that the coordinate is deleted.

\proclaim{Lemma 0 (Definition and properties of the $\SS_k$ [Bai])}
For each $k\in \{0, \ldots, n-1\}$ let
$\SS_k$   be the convolution operator
$$
\SS_k \phi(x) =(-1)^ n
\int_y \sigma_k (x,y) \wedge \phi(y) \, ,    \tag{2.7}
$$
where $\sigma_k (x,y)$ is a $k$-form in $x\in \real^n$
and an $n-k-1$-form in $y\in \real^n$,   measurable on $\real^{2n}$,
defined by  the decomposition ($\sigma(x-y)$ denotes $(F^* \sigma)(x,y)$
where $F(x,y)=x-y$)
$$
\sigma(x-y)=\sum_{k=0}^{n-1} (-1)^{n} \sigma_k(x,y) \, .  \tag{2.8}
$$

\noindent (1) 
For every bounded open ball $D \subset \real^n$, the form
$\sigma_k(\cdot, y)$ belongs to $\AA_{k,L^s(D)}$
for each $1 \le s < n/(n-1)$ and
$
\sup_{y \in D} \|\sigma_k(\cdot, y)\|_{L^s(D)} < \infty$ 
for all  $k$.
Outside of the diagonal $x=y$ in $\real^{2n}$,
$\sigma_k(x,y)$ has $C^\infty$ coefficients.

\noindent (2)
The following homotopy equation holds:
$
d_{k-1} \SS_{k-1} + \SS_k d_k = \Id\quad
\text{on } \AA^C_{k,C^2} 
$.
\endproclaim

\demo{Proof of Lemma 0}
By definition, the singular  kernel $\sigma_k(x,y)$
of $\SS_k$ can be written $\sigma_k(x,y) =
\sum_{\vec \imath} \tau_{k, \vec \imath} (x-y) dy^{\vec \imath} $, where
the coefficients of each $k$-form
$\tau_{k,\vec \imath}(u)$ are linear combinations of:
$$
{u_j \over \|u\|^n} \, , j=1, \ldots, n \, , \tag{2.9}
$$
where $\|u\|$ denotes the Euclidean norm of $u\in \real^n$.
The fact that $ \sigma_k(x,y)$ is $C^\infty$ outside of the
diagonal is obvious.

It was pointed out in [Bai] that the
singularity in \thetag{2.9} belongs to $L^{s}(D)$ for each
$1\le s < n/(n-1)$ and every bounded ball $D$.
The integrability claim in Lemma 0 easily  follows from this.
Claim \therosteritem{2} is not difficult (see [Bai]).
\qed
\enddemo

Finally, the kneading operators are defined (in the sense
of formal power series with coefficients linear operators)
as
$$
\DD_k(z)= z \NN_k (\Id-z \MM_k)^{-1} \SS_k
=\sum_{\ell=0}^\infty z^{\ell+1} \NN_k \MM_k^\ell \SS_k\, ,
\quad k=0, \ldots, n-1 \, .
\tag{2.10}
$$

\smallskip
{\bf The higher dimensional Milnor-Thurston  formula}
\smallskip

In order to state the results from [Bai], let us view the kneading
operators as kernel operators.
By Lemma~0, since both $\MM_k$ and $\NN_k$ act  boundedly on
$\AA_{k, L^q}$, we may  write
$$
\DD_k(z) \phi (x)=(-1)^n
\sum_{\ell=0}^\infty  z^{\ell+1} \int_y (\NN_k  \MM_k^\ell)_x
 \sigma_k(x,y) \wedge \phi(y) \, .\tag{2.11}
$$
(The subscript $x$ means that the operator acts on the
$x$ variable.)
Hence, $\DD_k(z)$ has a formal power series kernel
$\sum_{\ell \ge 0} z^{\ell+1} \KK_{k, \ell}(x,y)$
(just set $\KK_{k, \ell}(x,y)= (\NN_k  \MM_k^\ell)_x
 \sigma_k(x,y)$).
The kernels of iterates $\DD^m_k(z)$ are obtained
in a similar way and analogous expressions
may be written for $\widehat \DD_k(z)$ and $\DD^*_k(z)$.

Baillif [Bai] proved the following lemma and theorem
(see also [Ba2, \S3]) using transversality and elementary properties
of the $\sigma_k(x,y)$:

\proclaim{Lemma 1 [Bai]}
For each $k=0, \ldots, n-1$
and $\ell \ge 0$, we have $\KK_{k,\ell}(x,x)\in L^1(\real^n)$,
so that we may define the 
formal flat trace of $\DD_k(z)$
by integrating the kernel along the diagonal
(with the correct sign):
$$
\tr^\flat \DD_k(z)=
\sum_{\ell \ge 0} z^{\ell+1} (-1)^{(k+1)}
\int\KK_{k,\ell}(x,x)   \, . \tag{2.12}
$$
The same integrability property holds for the kernels of
$\DD_k(z)^m$ for all $m\ge 1$.
Therefore,  the traces $\tr^\flat \DD_k^m(z)$
are power series with complex coefficients.   
\endproclaim

For $k=0, \ldots, n-1$, we may thus define a (formal) flat determinant
by
$$
\Det^\flat(\Id+\DD_k(z))=
\exp -\sum_{m=1}^\infty {z^m \over m} \tr^\flat \DD_k(z)^m \, . 
$$

The coefficients   of $\tr^\flat \DD_k(z)$
coincide with the Atiyah-Bott flat traces of the
corresponding kernel operators if the $\psi_\omega$ and $g_\omega$
are $C^\infty$, see [Bai].

\proclaim{Theorem 2 (Milnor-Thurston identity) [Bai]}
In the sense of formal power series
$$
\eqalign{
\Det^\#(\Id-z \MM)
&=\prod_{k=0}^{n-1} \Det^\flat(\Id+\DD_k(z))^{(-1)^{k+1}}\, .
}\tag{2.13}
$$
\endproclaim

We shall not need the following result,
although we shall exploit and revisit key ideas in its proof (see [Bai, Lemma~6.2])
in Lemma~6 below:

\proclaim{Theorem  [Bai]}
There is $\delta > 0$ so that in the disc $\{ |z| < \delta \}$,
for all $m  \ge n/2+1$ 
the kernel of $\DD^m_k(z)$  has
coefficients in $L^2(dx\times dy)$.
In particular,  the regularised determinant of order
$n+1$ (if $n$ is odd) or $n+2$ (if $n$ is even)
of $\DD_k(z)$ on $L^2$ is holomorphic
in the disc of radius $\delta$.
\endproclaim

\smallskip
{\bf Harmonic analysis and algebra with the $\SS_k$}
\smallskip

\proclaim{Lemma 3 (Properties of the $\SS_k$ and $d_k\SS_k$: analysis)}
Let
$k\in \{0, \ldots, n-1\}$.

\noindent (1) If $1 < q < n$ then
$\SS_k$ is  bounded  from
$\AA_{k+1, L^q}$ to $\AA_{k, L^{q'}}$ for $q'=qn/(n-q)$.
For $1\le q\le n$ and any continuous
compactly supported function $\chi$, the
operator $\chi \SS_k \chi$ is bounded  from
$\AA_{k+1, L^q}$ to $\AA_{k, L^{q''}}$ for all $q<q''< qn/(n-q)$. 

\noindent (2)
$d_k \SS_k$ extends to a bounded operator
on $\AA_{k+1, L^q}$ for $1 < q < \infty$. In fact,
for each $i$, the operator $\partial_i \SS_k$ extends
to a bounded operator from $\AA_{k+1, L^q}$
to $\AA_{k, L^q}$.
\endproclaim

Note that $q' > q>1$, in fact $q' > \eta q$ for each $1\le \eta < 1+q/(n-q)$.

\demo{Proof of Lemma 3}
The second claim of \therosteritem{1} follows from
usual properties of the convolution (see e.g. [Sch, p.151])
since $\chi \sigma_k \in L^t$ for all $1\le t < n/(n-1)$.

Noting that the kernel $\sigma_k$ is
of weak-type ${n\over n-1}$, we may apply the
Hardy-Littlewood-Sobolev fractional integration theorem
(see e.g. [St, V.Theorem~ 1 and Comment~1.4];
the comment explains why Theorem~ 1 extends from the
Riesz potential to more general weak-type kernels),
which yields that $\SS_k$ is a bounded operator from
$\AA_{k+1,L^q}$ to $\AA_{k,L^{q'}}$ if $1 < q < n$.

It remains to prove \therosteritem{2}, i.e.,  study $\partial_i\SS_k$.
We first claim that $\partial_i \SS_k$ can be written as
the sum of, on the one hand a convolution
operator with kernel a  form having coefficients
which are linear combinations of expressions of the type
$$
{x_i x_j \over \|x\|^{n+2}} \, , {\|x\|^2-n x_i^2\over \|x\|^{n+2}} \, ,
1\le i,j\le n \, , \tag{2.14}
$$
and on the other hand
a distribution $\nu_k$  which extends to a bounded
operator on  $\AA_{k+1,L^q}$.
Indeed, computing the functional partial derivatives
$\partial_i={\partial\over \partial x_i}$ of \thetag{2.9} produces
\thetag{2.14},
so that we only need to study the distributional
contribution to $\partial_i \SS_k$. For this,
let us take a kernel singularity
$x_j/\|x\|^{n}$, and  for each $t>0$
let $\chi_t:\real\to \real$ be a $C^\infty$ function,
identically zero on $[-t,t]$ and which is identically
equal to $1$ outside of $[-2t,2t]$. We may assume that
$\sup |\chi'|\le 2/t$. Let us consider
$$
\partial_i
\left ( \chi_t(\|x\|) x_j \over \|x\|^{n} \right ) =
{x_j \partial_i (\chi_t(\|x\|)) \over \|x\|^{n}} +
\chi_t(\|x\|) \partial_i \left ({x_j\over  \|x\|^{n}} \right ) \, .
$$
Since the second term of the above sum converges to
the already mentioned functional partial derivative as $t\to 0$, we should
check that
the distribution with kernel corresponding to the first term,
which is just $\chi'_t(\|x\|) x_i x_j/\|x\|^{n+1}$,
acts boundedly on $\AA_{k+1,L^q}$, uniformly in
$t\ge 0$. By density and completeness, it is enough to
consider $\phi \in \AA_{k+1,L^q}\cap C^\infty$
and take $\varphi=\phi_\ell$ for $\ell\in I(k+1)$.
We may formally write for every $w\in \real^n$
$$
\eqalign
{
|\nu_k(\varphi)(w)| &=
\left |
\int_{\real^n} \chi'_t(\|x\|) {x_i x_j \over \|x\|^{n+1}}\varphi(w-x) \, dx
\right |\le
{2\over t} \int_{\|x\|\le t}
\left | {x_i x_j\over \|x\|^{n+1}}\varphi(w-x)\, dx
\right |\cr
&\le
{2\over t} \sup_{\|x\| \le t} |\varphi(w+x)|
\int_{\|x\|\le t} \biggl | {x_i x_j\over \|x\|^{n+1}}\biggr | \, dx    \, .
}\tag{2.15}
$$
Going to polar coordinates, it is not difficult
to see that  there is a constant $C>0$
so that $\int_{\| x\|\le t} |x_i x_j| \|x\|^{-n-1}
\, dx \le C t$. We end our analysis of the
distributional contribution by noting that
$
\lim_{t\to 0} \int \sup_{\|x\|\le t} |\varphi(w+x)|^q \, dw
=\int |\varphi(w)|^q \, dw 
$.

(Since our assumptions imply that $|\varphi(w+x)-\varphi(w)|
\le t \sup \|D \varphi\|$, the above is easily checked.)
As pointed out to us by S.~Gou\"ezel, a more careful
analysis shows that in fact $\nu_k$ is zero if
$i\ne j$ and is a scalar multiple of the Dirac mass if $i=j$.

We next observe that the
expressions \thetag{2.14} in the functional
term of the kernel of $\partial_j \SS_k$  exhibit
the same kind of singularity as the Riesz
transform. More precisely, this kernel is a linear
combination of $\Omega_\ell(x)\over \|x\|^n$ where
each $\Omega_\ell$
is homogeneous of degree zero, has vanishing integral
on the unit $(n-1)$-dimensional sphere  and is $C^1$
on this sphere. We may thus apply  e.g.  Theorem ~3
in Chapter II of [St], which immediately guarantees
that  $d_k \SS_k$ extends to a bounded operator
on $\AA_{k+1,L^q}$, for $1< q < \infty$. 
\qed\enddemo

\proclaim{Lemma 4 (Properties of the $\SS_k$: algebra)}
Let $1< q < \infty$ .

\noindent (1) 
$\SS_{k-2}  \SS_{k-1}=0$ on  $\AA_{k,L^q}$.

\noindent (2)
On $\AA_{k, L^q}$, and suppressing the indices for simplicity,
$d \SS=(d \SS)^2$, and $( \SS d)^2= \SS d$,
$d \SS  \SS d=0= \SS d d\SS$. In other words, $d\SS$ and $\SS d$ are two
orthogonal bounded projectors onto $\Imm d=\Ker d$,
respectively  $\Imm  \SS=\Ker \SS$, in $\AA_{k, L^q}=\Imm d \oplus \Ker \SS
=\Imm \SS \oplus \Ker d$.
\endproclaim

\demo{Proof of Lemma 4}
The facts that $\SS d$ is onto
$\Ker \SS$ and
that $\Ker  \SS \subset \Imm  \SS$ do
not depend on $ \SS^2=0$:
Use $\phi=\SS d\phi$ if $\phi \in \Ker  \SS$.

To prove that $\SS^2=0$,
we shall give an equivalent definition of $\SS$, which was indicated to us
by D.~ Ruelle.
Using Lemmas ~0 and ~3, it suffices to show the result on $C^\infty$ compactly supported
forms (by density).
Recall that there is a scalar (or hermitian) product on
the space of compactly supported $C^\infty$
$k$-forms $\AA^C_{k,C^\infty}$, defined by:
$$
 < \phi(x)dx_{i_1}\wedge\cdots\wedge dx_{i_k} \, |\,
\varphi(x)dx_{j_1}\wedge\cdots\wedge dx_{j_k}> =
 \left\{
 \eqalign{
 &\int \phi\cdot\varphi  \, \text{ if $i_1,\dots,i_k=j_1,\dots,j_k$} \cr
 &0 \, \text{ otherwise.}
 }
 \right. \tag{2.16}
$$
(The scalar product is  still defined if only one of $\phi,\varphi$ is compactly supported.)
Now, $d:A^C_{k,C^\infty}\to \AA_{k+1,C^\infty}$ defines a dual
$d^*:\AA^C_{k+1,C^\infty}\to \AA_{k,C^\infty}$ by
$
<d^*\phi,\varphi>=<\phi,d\varphi>\, .
$ 

The Laplace-Beltrami operator is defined by $\Delta = dd^*+d^*d$ where $d$ is
the exterior derivative of forms.
In $\real^n$, we have (see, e.g., [Li]):
$$
\Delta \left[ \sum_{j\in I(k)}\phi_j(x) dx_j \right] = -\sum_{j\in I(k)}
(\sum_{i=1}^n \frac{\partial^2}{\partial x_i^2}
\phi_j)(x) dx_j    \, .\tag{2.17}
$$
Define the Green kernel to be
$E(x)=-\frac{\Gamma(n/2)}{(n-2)2(\pi)^{n/2}}\frac{1}{||x||^{n-2}}dx_1
\wedge\cdots\wedge dx_n$
if $n\ge 3$, and
$E(x)=\frac{1}{2\pi}\log(||x||)dx_1\wedge dx_2$ if $n=2$.
Then, let $E_k(x,y)$ be the $k$-form in $x$ and $n-k$-form in $y$ such that
$E(x-y)=\sum_{k=0}^n (-1)^{n}E_k(x,y)$. It is well known
(see [Sch]) and has been used in [Bai] that
$\Delta E = \delta(x) dx_1\wedge\cdots\wedge dx_n$ (as a current acting on
$C^\infty$
forms vanishing
at infinity, or  compactly supported  $C^\infty$forms). One thus sees
easily that the operator $G_k:\phi(x)\mapsto\int_y E_k(x,y)\wedge\phi(y)$
satisfies
$\Delta G_k=id$ on $\AA^C_{k,C^\infty}$.
Using that $\frac{\partial^2}{\partial y_i^2}E(x-y)=\frac{\partial^2}
{\partial x_i^2}E(x-y)$,
one  checks that $G_k\Delta=id$ on the same space.

The operator $\SS_k$ defined in [Bai] is in fact $d^*G_{k+1}$. Indeed,
$$
d^*\left(\phi(x)dx_1\wedge\cdots\wedge dx_n\right)=
\sum_{j=1}^k(-1)^{j+1}\partial_j\phi(x)dx_1\wedge\cdots{d\hat{x}_j}
\cdots\wedge dx_n \, . \tag{2.18}
$$
Thus, $d^*E(x)=\sigma(x)$, where $\sigma$ was defined in \thetag{2.6}.
Since $\SS_k$ is defined in  \thetag{2.7}
as the convolution with $\sigma_k$, the equality is
immediate.

By Lemma 3, the composition $\SS_{k-1}\SS_{k}=d^*G_kd^*G_{k+1}$ is well
defined on $L^q$ forms, and thus on $\AA^C_{k,C^\infty}$.
We shall show that $d^*G_k=G_{k-1}d^*$, which implies $\SS_{k-1}\circ\SS_{k}=0$
since $d^*d^*=0$.
Since $E(x)$ is an $n$-form, $d E(x)=0$. Using $d=d^x+d^y$, we obtain
$d^xE_{k-1}(x,y)=-d^yE_k(x,y)$. Integrating by parts,
$dG_{k-1}=G_{k}d$
(recall that
$E_k(x,y)$ is an $n$-form). Notice now that $\Delta$ is auto-dual by
definition,
and thus $G_k$ is also auto-dual.
Hence, $d^*G_k=G_{k-1}d^*$.
\qed\enddemo

\smallskip

\head
3. Spectral interpretation of the zeroes of the sharp determinant
\endhead

{\bf The Axioms}

\smallskip

Let $r$, $\Omega$,
$\psi_\omega$, $g_\omega$, $K=\{ x\in \real^n \mid |x|\le T\}$, $\MM_k$,
$\SS_k$, $\sigma_k$, $\NN_k$, and
$\DD_k(z)$ be the objects from   Section~ 2.
Let $1_K$ be the characteristic function of $K$.
We fix once and for all
$$
K' =\{ x\in \real^n \mid |x|< 2 T\}\, ,
$$
and a radial
$C^\infty$ function
$\chi_K=\chi_{K,K'}$  supported in $K'$ and identically equal
to $1$ on $K$.  Note that $\NN_k$ and $\MM_k$
are bounded from $k$-forms to $k+1$-forms,
respectively $k$-forms, with coefficients in $L^q(K')$. 
We use the notation $\rho(\PP)$, $\rho_{ess} (\PP)$ for
the spectral and essential spectral radii of a bounded linear
operator $\PP$.
We are now ready to state our three  assumptions:

\smallskip

{\bf Axiom 1:}
For $k=0, \ldots, n$ and $1 < t < \infty$,
there are Banach spaces $\BB_{k,t}$ of distributions  on $\real^n$
so that:
\roster
\item
There are  real self-adjoint
invertible pseudodifferential
convolution operators $\widetilde \JJ_{k}$, $\widetilde \JJ_{k}^{-1}$ of order 
(at most) $r$
so that   $ \BB_{k, t}$
is defined by the Banach space isomorphism
$$
\widetilde \JJ_{k}(L^t(\real^n)) = \BB_{k, t} 
\, .\tag{3.1}
$$
\item
If $\psi$ is a $C^r$ local diffeomorphism
and and $g$ is a $C^{r}$ function supported in $K$
then   $\varphi\mapsto g \cdot \psi ^* (\varphi)$ is a bounded
operator on each $\BB_{k,t}$.
\item
Each operator $d_k :\AA_{k, \BB_{k,t}}
\to \AA_{k+1, \BB_{k+1,t}}$ is bounded.
\endroster

{\bf Axiom 2:} There is  $0\le \widetilde R < 1$ so that for $k=0, \ldots, n$:
\roster
\item
 $\MM_k$ is bounded on $\AA_{k,\BB_{k,t}}$ 
with $\rho (\MM_k|_{\AA_{\BB_{k,t}}}) \le 1$,
$\rho_{ess} (\MM_k|_{\AA_{\BB_{k,t}}}) < \widetilde R$, 
and  there are no eigenvalues of
modulus $\widetilde R$, for each $1< t< \infty$.
\item
Let $\Pi_{k,t}$ be the (finite-rank) spectral projector on $\AA_{k,\BB_{k,t}}$
associated to the spectrum of $\MM_k$ outside the disc of
radius $\widetilde R$. Let $\{\varphi_{k,t,s}\}$ be a basis for the corresponding
generalised eigenspace. (Our assumptions imply that the dimension
does not depend on $1 <t<\infty$.)
Then, there is $C>0$ so that
$\max_{k,s} \sup_{1 < t < \infty}
 \| \varphi_{k,t,s} \|_{\BB_{k,t}}\le C$, and, for all $k$,  $j$ , we have
$$
 \sup_{1 < t < \infty}
 \| \MM_k^j \|_{\BB_{k,t}}\le C\, ,
\quad  \sup_{1< t < \infty} \|\MM_k ^j (\Id-\Pi_{k,t})\|_{\BB_{k,t}}
 \le  C \widetilde R^j \, .
$$
\endroster

{\bf Axiom 3:} 
There is $\EE_0 > 0$ so that,  for 
all $0<\EE < \EE_0$, each $j\ge 1$
and each admissible composition $\psi^{j}_{\vec \omega}$,  letting
$V_{\vec \omega,\EE}$ be the set of $x$ 
in the domain of definition of  $\psi^{j}_{\vec \omega}$
so that
$\|\psi^{j+\ell}_{\vec \omega}(x)-x\| < \EE$, the map
$\psi_{\vec \omega}^{j+\ell}- \Id$ is injective on $V_{\vec \omega,\EE}$.

Also, for each
$\eta >0$ there is $C$ so that for all $k$, every sequence
$\{\EE_j\}$ with $0<\EE_j <\EE_0$, and each $k$-form 
$\bold 1_k= dx_{i_1} \wedge \cdots \wedge dx_{i_k}$
$$
\sum_{\vec \omega \in \Omega^j}
\sup_{x \in V_{\vec \omega,\EE_j}}
{|g_{\vec \omega}^{(j)} (\psi^{j}_{\vec \omega} x)| \over |\det (D\psi_{\vec \omega}^{j}(x) -\Id ) | }
\bigl | (\psi_{\vec \omega}^{j})^* (\bold 1_k) (x) / \bold 1_k(x) \bigr | 
\le C \exp(\eta j) \, , \, \forall j\, .
$$


\smallskip

\remark{Consequences of Axiom 1}
$C^{\infty}$ functions with compact support
are contained in $\BB_{k,t}$ and 
distributions in $\BB_{k,t}$ have order at most $r$.
Also,  $\BB_{k,t}\subset \BB_{k,t'}$
boundedly if $t >t'$.

Since $\widetilde \JJ_k$ is a convolution 
operator it commutes  (at most up to  sign)
with $\SS$, $d\SS$,  $\SS d$ and with convolution operators
having a smooth kernel $\delta_\epsilon(\|x-y\|)$.
Thus, by Lemma ~3,  for all $1< t < n$ the operators
$\SS_k :\AA_{k+1, \BB_{k,t}}\to \AA_{k, \BB_{k,tn/(n-t)}}$ are bounded.
Also, for all $1 < t <\infty$  the operators
 $\SS_k:\AA_{k+1, \BB_{k+1,t}}\to \AA_{k, \BB_{k,t}}$ are bounded,
and   both operators $d_k  \SS_k$, $\SS_{k+1}d_{k+1}$ are
bounded on $\AA_{k+1, \BB_{k+1,t}}$.
Finally, there are $C^\infty$ mollifiers $\delta_\epsilon(\|x-y\|)$
on $K'\times K'$ 
with $\|\int \delta_\epsilon (\cdot,y)\varphi(y) \, dy - \varphi\|_{\BB_{k,t}}\to 0$
as $\epsilon \to 0$ for all $\varphi$ [Ad, Lemma 2.18].

By  Lemma~0 \therosteritem{1}, for  each fixed $y\in K$,
the form $\sigma^{(r)}_k(x,y)=\widetilde \JJ_{k}
(1_K(\cdot) \sigma_k(\cdot, y))$ belongs to $\AA_{k,\BB_{k,t}}$
for all $1 \le t < n/(n-1)$ and each $k$, with
$$
\max_k \sup_{y\in K}   \|  \widetilde \JJ_{k}
(1_K(\cdot) \sigma_k(\cdot, y)) \|_{\AA_{k,\BB_{k,t}}} < \infty  \, ,
\forall  \, 1 \le t < n/(n-1)\, .
\tag{3.2}
$$
Write $\SS^{(r)}_k$ for the  convolution operator
corresponding to $\sigma^{(r)}_k(x,y)$. Up to 
a sign, this is just $1_K \widetilde K \JJ_k \SS_k$. 
Hence, there is $\upsilon_k\in\{+1,-1\}$ so that
for all $1<t<\infty$
$$
1_K \SS_k \widetilde \JJ_{k}=\upsilon_k  \SS^{(r)}_k\, ,\quad 
\hbox{ on } L^t(\real^n) \, .
\tag{3.3}
$$
Since $\widetilde \JJ^{-1}$ is real and self-adjoint,
for all  compactly supported
$\varphi \in L^\infty(\real^n)$, all 
$\psi \in \BB_{k,t}$
$$
\int (\widetilde \JJ^{-1}_k \varphi)(x) \psi  (x) dx
=\int \varphi (x) (\widetilde \JJ^{-1}_k \psi)(x) dx \, . \tag{3.4}
$$
\endremark

\smallskip


{\bf Preliminary Step: The $L$-regularised kneading determinants $\DD_{k,L}(z)$}
\smallskip

For each $\ell \ge1$
the $\ell$-regularised version of $\MM_k$ is defined to be the following
convergent power series
with operator (on $\AA_{k,\BB_{k,t}}$, e.g.) coefficients:
$$
\MM_{k,\ell}(z) =\Id - (\Id - z\MM_k)
\exp\biggl (\sum_{j=1}^{\ell-1}{z^j\MM_k^j\over j} \biggr ) \, . \tag{3.6}
$$
Note that $\MM_{k,1}(z)=z\MM_k$ and   $\MM_{k,\ell}(z)$ converges for
all $z$ and all $\ell$.
We shall also use
$
\NN_{k,\ell}(z)= d\MM_{k,\ell}(z)- \MM_{k+1,\ell}(z) d
$.
(Note that $\NN_{k,1}(z)=z\NN_k$.)

\proclaim{Lemma 5} Assume Axioms 1 and 2.
For each $C > 1$ and $\xi < 1/\widetilde R$ there is $L \ge 1$ so that for each
$|z| < \xi$ and all $t$, $k$, we have
$\|\MM_{k,L} (z)|_{Im(\Id-\Pi_{k,t})}\|_{\AA_{k,\BB_{k,t}}}\le
{|z|^L  \widetilde R^L\over C}$ and
$$
\|\NN_{k,L} (z)|_{Ker(\Pi_{k,t})\cap Ker(\Pi_{k+1,t} \circ d_k)}\|_{\AA_{k,\BB_{k,t}}, 
\AA_{k+1,\BB_{k+1,t}}}\le {|z|^L  \widetilde R^L\over C} \, .
$$
Up to taking a slightly larger value of $\widetilde R$,
we may assume that $L$ does not depend on $\xi$.
\endproclaim

\demo{Proof of Lemma 5}
Formally,
$$
\MM_{k,\ell}(z)= \Id -(\Id-z\MM_k)\exp[-\log(\Id-z \MM_k)-\sum_{j\ge \ell}
{z^j\MM_k^j\over j}]
=\Id-\exp\sum_{j\ge \ell}- {z^j \MM_k^j\over j} \, .
$$

Then, we use that for any $0<\theta < 1$ and $C > 1$
there is $L$ so that for all $\ell\ge L$
$$
\sum_{j\ge \ell} {\theta^j\over j} \le { \theta^\ell\over \ell}
\sum_{j\ge 0} \theta^j\le {\theta^\ell \over \ell (1-\theta)}
\le {\theta^\ell\over 2 C} \, ,
\hbox{ while }
1-\exp(-{\theta^\ell\over 2C}) \le {\theta^\ell\over C} \, .
$$
This gives the bound for $\MM_{k,\ell}(z)$ by Axiom~2\therosteritem{1}.
To estimate $\NN_{k, \ell}(z)$ 
note that
$$
\NN_{k,\ell}(z)
=\Id-\exp\sum_{j\ge \ell}- {z^j \NN_k^{(j)}\over j} \, ,\tag{3.8}
$$
with $\NN^{(j)}_k=d_k\MM_k^j - \MM^j_{k+1} d_k$,
and use Axiom~1\therosteritem{3}.
\qed
\enddemo

Whenever the value of $t$ is clear
from the context, we write $\Pi_{k}$ instead of $\Pi_{k,t}$ for simplicity.
Note that
$\Pi_k =\Pi_k \chi_K=\chi_K \Pi_k$.
It follows from Lemma~5 that the essential spectral
radius of $\MM_{k,\ell}(z)$ on any $\AA_{k,\BB_{k,t}}$ is not larger
than $|z|^\ell \widetilde R^\ell< 1$ for all $\ell$ and every $|z|< 1/ \widetilde R$.
It is then an easy algebraic exercise to see that for any $|z|< 1/ \widetilde R$,
the complex number $1/z$ is an eigenvalue of $\MM_k$ (on $\AA_{k,\BB_{k,t}}$)
of algebraic multiplicity $m$ if and only $1$ is an eigenvalue of
multiplicity $m$ for $\MM_{k,\ell}(z)$ (on $\AA_{k,\BB_{k,t}}$), and
in particular, for all $\ell$, setting $R=\widetilde R^{1/n}$,
$$
V_{k,t}=\{ |z| < 1/ R \mid 1/z \notin \sp \MM_k \}
= \{ |z| < 1/  R \mid 1 \notin \sp \MM_{k,\ell}(z) \}\, .\tag{3.7}
$$
(Just use that $\MM_k$ commutes with each $\exp(z^j \MM_k^j/j)$.)
One can use the same basis of generalised eigenvectors for
both eigenspaces, and $\Pi_k$ commutes
with each $\MM_{k,\ell}$. (See e.g. [GGK] for analogous results in the case
when a power of $\MM_k$ is Hilbert-Schmidt.)

We shall work with the regularised kneading operators
(for suitably large $\ell$)
$$
\DD_{k,\ell}(z)= \NN_{k,\ell}(z) (\Id- \MM_{k,\ell}(z))^{-1} \SS_k \, .
\tag{3.9}
$$
We explain next why Lemma 1  and a modified
version of Theorem ~2 hold for the  $\DD_{k,\ell}(z)$. 
The replacement of $z\MM_k$  and $z\NN_k$ by $\NN_{k,\ell}(z)$,
$\MM_{k,\ell}(z)$ does not cause any problems in
Baillif's [Bai, Ba2] proof of Lemma~ 1, 
since $\MM_{k,\ell}(z)$
is just an entire series with coefficient transfer operators 
(acting on $\AA_{k,L^q(K')}$).
However, since  \thetag{2.12} implies
$$
\Det^\flat(\Id-\MM_{k,\ell}(z))=\Det^\flat(\Id-z\MM_k)
\exp\biggl (\sum_{j=1}^{\ell-1}{z^j\over j} \tr^\flat \MM_k^j\biggr )\, ,
$$
the first equality in Theorem~2 must be replaced by
$$
\Det^\#(\Id-z \MM)\, \exp\biggl 
( \sum_{j=1}^{\ell-1} {z^j\over j} \tr^\#\MM^j\biggr )
=\prod_{k=0}^{n-1} \Det^\flat(\Id+\DD_{k,\ell}(z))^{(-1)^{k+1}}
\, .\tag{3.10}
$$
The additional factor 
$\exp\sum_{j=1}^{\ell-1} {z^j\over j} \tr^\#\MM^j$ is
clearly an entire and non vanishing function of $z$ for each $\ell$.
The final useful property of the regularised kneading determinants is:

\proclaim{Corollary of Lemma 5}
Assume Axioms 1 and 2. 
For each
$z\in V_{k,t}$ (see \thetag{3.7})
and all $0\le k\le n-1$,    $1 < t < \infty$,
the essential spectral radius of
$\DD_{k,L}(z)$ on $\AA_{k+1,\BB_{k+1,t}}$ 
goes to zero exponentially fast as $L\to \infty$, in particular it is $<1$ for large~ 
$L$.
\endproclaim

\demo{Proof of the Corollary}
Let $\widetilde \Pi_k$ be a projector onto
the finite-dimensional space 
$\{ \psi \in (\Id-\Pi_k)(\AA_{k,\BB_{k,t}})  \mid d_k\psi \in \hbox{Im}\, \Pi_{k+1}\}$.
Use $\Pi_k\MM_{k,L}(z)= \MM_{k,L}(z)\Pi_k$ to get the decomposition
$$
\eqalign{
&\DD_{k,L}(z)=\NN_{k,L}(z) (\Id-\widetilde \Pi_k)(\Id- \MM_{k,L}(z))^{-1} (\Id-\Pi_k) \SS_k\cr
&\quad+\NN_{k,L}(z) \widetilde \Pi_k(\Id- \MM_{k,L}(z))^{-1} (\Id-\Pi_k) \SS_k
+\NN_{k,L}(z) (\Id- \MM_{k,L}(z))^{-1}\Pi_k \SS_k 
}
$$
into an operator of arbitrarily small spectral radius 
(by  $|z|< 1/R$ and Lemma~5)
and a finite rank operator (since $z\in V_{k,t}$). 
\qed
\enddemo

\smallskip

{\bf Operators $\DD^{(r)}_{k,L}(z)$ }

\smallskip

We next introduce auxiliary operators $\DD_{k,L}^{(r)}(z)$ on
$\AA_{k+1, L^2(K')}$, extending a one-dimensional construction of
Ruelle [Ru4]. Their iterates will be trace-class, and their
traces  will coincide with the formal flat traces of iterates
of the $\DD_{k,L}(z)$
in the sense of power series. This will allow us to prove
the following crucial lemma:

\proclaim{Lemma 6 (Meromorphic extension of $\Det^\flat(\Id+ \DD_{k,L}(z))$)}
Assume Axioms 1--2 and~3. Set
$\BB_k=\BB_{k,2}$ for  $k=0, \ldots, n$.
For $k=0, \ldots, n-1$
and all large enough ~ $L$:

\noindent (1)
$\Det^\flat(\Id+ \DD_{k,L}(z))$ extends holomorphically to
$$V_k=\{ |z| < 1/R \mid 1/z \notin \sp ({\MM_k}|_{\AA_{k,\BB_{k}}}) \}
=\{ |z| < 1/  R \mid 1 \notin \sp (\MM_{k,\ell}(z)|_{\AA_{k,\BB_{k}}}) \}\, .$$

\noindent (2)
If  $|z| < 1/R$ and $1/z \in \sp({\MM_k}|_{\AA_{k,\BB_{k }}})$, then
$\Det^\flat(\Id+\DD_{k,L}(z))$ is meromorphic at $z$ with
a pole of order at most the algebraic multiplicity
of the eigenvalue.

\endproclaim

\demo{Proof of Lemma 6}
Recall  $\sigma^{(r)}_k(x,y)$ and $\SS^{(r)}_k$ from Axiom~1,
and that $\sigma^{(r)}_k(\cdot,y) \in \AA_{k,\BB_{k,t}}$,
for $1 \le t < n/(n-1)$ and $y\in K$.
Set:
$$
\eqalign{
\DD_{k,L}^{(r)}(z)=\upsilon_k
 (\widetilde \JJ_{k})^{-1} \NN_{k,L}(z) (\Id-\MM_{k,L}(z))^{-1}  \SS^{(r)}_k\, .
}
$$

Axioms 1 and 2 first imply that $z \mapsto \DD_{k,L}^{(r)}(z)$ is a well-defined
map (taking values in the space of bounded linear
operators on $\AA_{k+1, L^2(K')}$), holomorphic in  $V_k$,
and meromorphic in the disc of radius $1/\widetilde R$, with possible poles
at the inverse eigenvalues of $\MM_k$ on $\AA_{k,\BB_{k,2 }}$
(the order of the
pole being at most the algebraic multiplicity of the eigenvalue).
We use here that, by Axiom~1(2-3), 
$\NN_{k,L}(z)$ is a $z$-entire function, bounded 
from $\AA_{k, \BB_{k,t}}$ to
$\AA_{k+1, \BB_{k+1,t}}$ for all $1<t<\infty$.

The kernel of $\DD_{k,L}^{(r)}(z)$ is
$$
\KK^\DD(x,y)=\upsilon_k (\widetilde \JJ_{x,k})^{-1} \NN_{x,k,L}(z)
 (\Id-\MM_{x,k,L}(z))^{-1} \sigma_k^{(r)} (x,y)  \tag{3.11}
$$
(where the transfer operators act on the $x$-variable).
The coefficients of
the kernel  $\KK^\DD(x,y)$ have no reason to be in $L^2(K\times K)$,
but  Axioms 1 and 2 imply that the coefficients of $\KK^\DD(\cdot, y)$  are
in $L^{t}(K')$ for all $1\le t < n/(n-1)$,
with the supremum over $y \in K$ of the $L^{t}$ norm bounded.

The proof of Lemma~6.2
in [Bai]  shows that for $z \in V_k$, and all $\ell > n/2$
the operator   $(\DD_{k,L}^{(r)}(z))^\ell$ 
is Hilbert-Schmidt on $\AA_{k+1, L^2(K')}$. In particular,
the regularised determinant of order $[n/2]+1$,
$$
\Det^{\hbox{reg}}_{[n/2]+1}(\Id+\DD_{k,L}^{(r)}(z))=
\exp -\sum_{\ell = [n/2]+1}^\infty {z^\ell\over \ell}
\tr^\flat (\DD_{k,L}^{(r)}(z))^\ell \, , 
$$  
is holomorphic on $V_k$.
Indeed,  $(\Id-\MM_{k, L}(z))^{-1}$ is holomorphic in
$V_k$ as a bounded operator on
$k$-forms with coefficients in $\BB_{k,t}$, all $t$,
and in particular for all values of $t$ between $t_0 < n/(n-1)$ and
$t_{[n/2]+1}\ge 2$ which appear along the successive iterations.

By \thetag{3.3}
$$
\DD^{(r)}_{k,L}(z)= (\widetilde \JJ_{k})^{-1}\NN_{k,L}(z)
(\Id-\MM_{k,L}(z))^{-1}\SS_k \widetilde \JJ_{k} \, . \tag{3.12}
$$
Recall that the wave front set of (the Schwartz kernel of) a pseudo-differential
operator is included in $\{(x,x,\xi,-\xi)\mid \xi \ne 0\}$ and that composition
with a pseudodifferential operator does not enlarge the wave-front set
(see e.g. [AG]).
Hence, using the assumptions on  $\widetilde \JJ_k$,  $(\widetilde \JJ_k)^{-1}$,
a convolution with a $C^\infty$ mollifier,
the transversality property \thetag{2.1} to restrict
to the diagonal, and Fubini (see Section 4 of [Bai]), 
the following equalities
between formal power series hold
$$
\Det^\flat (\Id+\DD_{k,L}^{(r)}(z))
=\Det^\flat (\Id+\DD_{k,L}(z)) \, ,
\, \, \forall \, k=0, \ldots n-1
\, .\tag{3.13}
$$

\smallskip
To show that the full flat determinant $\Det^\flat(\Id+\DD_{k,L}(z))$ 
is holomorphic in $V_k$, we shall use Axiom 3 to  see that
the power series for $\tr^\flat (\DD_{k,L}(z))^\ell$
for each $1\le \ell \le [n/2]$ is
holomorphic on $V_k$ (its exponential is thus holomorphic
and nonvanishing).
We consider the (hardest) case, i.e., $\ell=1$, leaving higher iterates
to the reader. Also, we
assume for simplicity that $k=0$ (dealing with $k$-forms only introduces
notational complications).  

\smallskip 

We first show that the flat trace of
$(\widetilde \JJ_0)^{-1}\NN_{0,L}(z)(\Id-\MM_{0,L}(z))^{-1} (\Id-\Pi_0) \SS^{(r)}_0$
is holomorphic in the disc of radius $1/R$.
Let $\TT$ denote the restriction to the diagonal $x=y$
in $\real^n \times \real ^n$.
For each $\EE < \EE_0$, let  $W_\EE$ be the set of 
$(x,y)\in \real^n \times \real ^n$ with
$\|x-y\| < \EE$ and
let $\chi^\EE(x-y)+(1-\chi^\EE(x-y))$ be a 
$C^\infty$ radial partition of unity of $\real^n\times \real^n $ subordinated 
to $W_\EE$ in the sense that
$\chi_\EE(u)$ depends only on
$\|u\|$, $\chi^\EE(u)=1$ if $\|u\| \le \EE$,
and $\chi^\EE(u)=0$ if $\|u\| \ge 2\EE$.
Let $\delta_\epsilon$ be $C^\infty$ mollifiers
(converging to the Dirac mass in $\real^n$
as $\epsilon\to 0$),  and
write $\sigma^{\epsilon,<\EE}$ for the convolution
of $\delta_\epsilon$ with $\chi^\EE\sigma_0$
and  $\sigma^{\epsilon,>\EE}$ for the convolution
of $\delta_\epsilon$ with $(1-\chi^\EE)\sigma_0$

Using the identity $(\Id-\MM_{k,L}(z))^{-1}
= \exp - [\sum_{i=0}^{L-1} z^i \MM^i_k ]\cdot (1- z\MM_k)^{-1}$,
it suffices to show that
 there is $C$ so that for all $j$ there is $\EE$
so that for each $\epsilon>0$, all $|z|< 1/R$, 
$$
\eqalign{
&\int_x | \TT [   (\delta_\epsilon \star(
\widetilde \JJ^{-1}_{0}\NN_{0,L}(z)
( \exp (- \sum_{i=0}^{L-1} z^i \MM^i_0)  \MM_{0}^j  (\Id-\Pi_0) 
\widetilde \JJ_{x,0}(\chi_K \sigma^{\epsilon,>\EE})(x,y))] | dx
\cr
&\qquad  +
\int_x \TT | [   (\delta_\epsilon \star(
\NN_{0,L}(z)
( \exp (- \sum_{i=0}^{L-1} z^i \MM^i_0)  \MM_{0}^j  (\Id-\Pi_0) 
\sigma^{\epsilon,<\EE})(x,y))]| dx\le
C R^j\, .}
$$

The entire series $\NN_{0,L}(z)
 \exp (- \sum_{i=0}^{L-1} z^i \MM^i_0)=\sum_{q\ge 1}z^q \QQ_q $ has
coefficients 
$$
\QQ_q = \sum_{q_1+q_2=q}
\kappa_{q_1, q_2} \NN_0^{(q_1)} \MM_0^{q_2}:
 \AA_{0,\BB_0} \to \AA_{1,\BB_{1}}\, ,
\quad \kappa_{q_1, q_2}\in \real \, ,
$$
so that for any $\widehat R$ there is $C$ so that
$|\kappa_{q_1, q_2}|< C \widehat R^{q_1+q_2}$.
We may write $\NN_0^{(q_1)} \MM_0^{q_2}\varphi=
\sum_{\vec \omega \in \Omega^{q_1+q_2}}
h_{q_1,q_2,\vec \omega} \wedge (\psi_{\vec \omega}^{q_1+q_2})^*\varphi$
where the $h_{q_1,q_2,\vec \omega}$ are one-forms with $C^{r-1}$
coefficients supported in the domains of the $\psi_{\vec \omega}$.

Write $\sigma^{\epsilon,<,>\EE}(x,y)=\sum_m \tau^{\epsilon,<,>\EE}_m(x-y)$.
Fix $j$ and consider  the (hypothetical) case $q_2=q$. 
The  modulus of  the contribution  of 
$ \MM_0^{q_2} \MM^j_0$ may be bounded by the sum over $m$ of
$$
\align
& \biggl | 
\int
\TT[ \delta_\epsilon \star \bigl (
\MM_{x,0}^{j+q} ( 
\tau_m^{\epsilon,<\EE})(x- y)  \bigr ) ]\, dx \biggr |
+
\biggl | 
\int
\TT[ \delta_\epsilon \star\bigl  (\MM_{x,0}^{j+q}
\Pi_{x,0} ( 
\tau_m^{\epsilon,<\EE})(x- y) \bigr )]\, dx \biggr |
\cr
&\quad + \biggl | 
\int \TT[ \delta_\epsilon \star\bigl (
\widetilde \JJ_0^{-1}\MM_{x,0}^{j+q}(\Id-\Pi_{x,0})
\widetilde \JJ_0   (1_K 
\tau_m^{\epsilon,>\EE})(x- y) \bigr ) ]\, dx  \biggr | \, .\tag{3.15}
\endalign 
$$
Recall that for $\EE < \EE_0$
as in Axiom~3,  $V_{\vec \omega,\EE}$ is the set of $x$ so that
$\|\psi^{j+q}_{\vec \omega}(x)-x\| < \EE$. 
By Axiom 3, writing $\Phi_{\vec \omega}$ for the inverse 
of $\psi_{\vec \omega}^{j +q}- \Id$ on $V_{\vec  \omega,\EE}$,
and performing the corresponding change of variable,
the first term in the above sum is bounded by
$$
\align
&
\biggl | \sum_{\vec \omega}
\int_{(\psi_{\vec \omega}^{ j+q}- \Id) V_{\vec \omega,\EE}}  
{  g^{(j+q)}_{\vec \omega}(u+\Phi_{\vec \omega}(u))   
\over  |\det (D\psi_{\vec \omega}^{j+q}(\Phi_{\vec \omega} (u)) -\Id ) |}
 (
\tau_m^{\epsilon,<\EE})(u) 
\,  du \biggr | \cr
&\quad\le
 \sum_{\vec \omega}
\sup_{V_{\vec \omega,\EE}}
{ |g^{(j+q)}_{\vec \omega}\circ \psi^{j+q}_{\vec \omega}   |
\over  |\det (D\psi_{\vec \omega}^{j+q}(x) -\Id ) |}
\cdot
\int  |\tau^{\epsilon,<\EE}_m(u)|\,  du 
\le C \exp(\eta (q+j))   \EE  \, .
\endalign
$$

By the consequences of Axiom~1, for all $x$, $y$ in $K$, we have
 $\widetilde \JJ_{0} ((1_K \tau_m^{\epsilon,<\EE})(\cdot- y) )(x)=
\tilde \upsilon_k \widetilde \JJ_{0}((1_K  \tau_m^{\epsilon,<\EE})(x- \cdot) )(y)$
for $\tilde \upsilon_k\in \{+1, -1\}$.
Hence, for the second term, use \thetag{3.4}
to get (in the case where all eigenvalues are semisimple)
$$
\align
&\biggl | 
\int \TT[\delta_\epsilon \star 
\sum_s  \lambda_s^{j+q} \varphi_s(x)
 \nu_s(
\widetilde \JJ_{y,0}^{-1} \widetilde \JJ_{0,\cdot} (1_K 
\tau_m^{\epsilon,<\EE})(\cdot- y) ) ]\, dx \biggr |
\cr
&\le
\biggl | 
\int \TT[\delta_\epsilon \star 
\sum_s   \lambda_s^{j+q} \varphi_s(x) \widetilde \JJ_{y,0}^{-1}
 (\nu_s(
 \widetilde \JJ_{0,\cdot}( 1_K 
\tau_m^{\epsilon,<\EE})(\cdot- y) ) )
]\, dx \biggr |\cr
&\le 
\sum_s 
\bigl | \int_K  \lambda_s^{j+q} \varphi_s(x) \widetilde \JJ_{x,0}^{-1}
 (\nu_s(\widetilde \JJ_{0,\cdot}(1_K  
\tau_m^{\epsilon,<\EE})(\cdot- x) ) )  \, dx \bigr |\cr
\allowdisplaybreak
&= 
\sum_s 
\bigl | \int_K   \lambda_s^{j+q}
\widetilde \JJ_{0}^{-1}( \varphi_s )
 \nu_s(\widetilde \JJ_{0,\cdot}( 1_K \chi^\EE\cdot
\tau_m^\epsilon)(\cdot- x)  )  \, dx \bigr |
\cr
\allowdisplaybreak
&\le C 
\sum_s  \sup| \widetilde \JJ_{0}^{-1}( \varphi_s )|\cdot
 \sup_x \|\widetilde \JJ_{0,\cdot}(1_K  
\tau_m^{\epsilon,<\EE})(\cdot- x) ) \|_{\BB_{0,t}}   \cr
&\le  C  \sup_{x\in K} \|( 
\tau_m^{\epsilon,<\EE})(\cdot- x)  \|_{L^t(K)}
 \le  C
\sup_x (\int  |\tau_m^{\epsilon,<\EE}(u-x)|^t \, du)^{1/t}
\le  C \EE
 \, ,
\endalign
$$
taking $t>1$ close to $1$ and using
$|\lambda_s|\le 1$.
Nilpotent contributions produce polynomial growth in $q+j$ which
gives a bound $C \exp(\eta (q+j)) \EE$ for arbitrarily small $\eta$.

If $a(x,y)$ is continuous then $\int_K |a(x,x)| \, dx 
\le \hbox{Vol}(K) \sup_{y}\sup_x |a(x,y)|$.
Recall also that if $\sup_{t\to \infty} \|b(\cdot,y)\|_{L^t(K)} < B(y)$, then 
$b(\cdot,y)$ is in $L^\infty(K)$
and $\sup_K|b(\cdot,y)|<B(y)$. 
The third term may thus be estimated by
(use both parts of Axiom 2)
$$
\align
&
C  \sup_{x,y \in K} | \delta_\epsilon \star(\widetilde \JJ_{x,0}^{-1}
\MM_{x,0}^{j+q}(\Id-\Pi_{x,0})\widetilde \JJ_{x,0}(
 \tau_m^{\epsilon,>\EE})(x- y))) |  \cr
&\le C  \sup_y \sup_{t\to \infty}
\| \MM_{0}^{j+q}(\Id-\Pi_{0})
\widetilde \JJ _0 (1_K\tau_m^{\epsilon,>\EE})(\cdot- y)\|_{\BB_{0,t}}\cr
\allowdisplaybreak
&\le C \widetilde R^{j +q} \sup_{y, t}
\bigl
 (\|\widetilde \JJ_0(1_K  \tau_m^{\epsilon,>\EE})(\cdot- y)\|_{\BB_{0,t}}
+ \| \Pi_0 \widetilde \JJ_0(1_K  \tau_m^{\epsilon,>\EE})(\cdot- y)\|_{\BB_{0,t}}
\bigr )\cr
\allowdisplaybreak
&\le C \widetilde R^{j +q}\sup_y \sup_{t\to \infty}
\biggl (\|\tau_m^{\epsilon,>\EE}(\cdot- y)\|_{L^t(K)}\cr
&\qquad \qquad \qquad\qquad+ \sum_s
|\nu_s[\widetilde \JJ_0(1_K  \tau_m^{\epsilon,>\EE})](\cdot- y)| \cdot
\| \widetilde \JJ_0^{-1} \varphi_s\|_{L^t(K)}
\biggr )
\cr  
&\le C \widetilde R^{j+q} 
\bigl (   \sup_{\|u\|\ge \EE} | \tau^{\epsilon,>\EE}_m(u)|  
+ C \sup_y \|1_K \tau_m^{\epsilon,>\EE}(\cdot- y)\|_{L^{t_0}}   \bigr )  \cr
&\le C\widetilde R ^{j+q} ( \EE^{-(n-1)} + C) \, , \cr
\endalign
$$
for some $1 < t_0 < n/(n-1)$.
Since $\|\NN^{(q_1)}_0\|\le C^{q_1}$,
 $\Omega$ is finite,  and
$ \|h_{q_1,q_2,\vec \omega}\|\le C^{q_1+q_2}$,  the case $q_1 \ge 1$ follows
from the bounds on $\kappa_{q_1,q_2}$.
Choosing $\EE=\widetilde R^{j/n}$  gives claim (1)
(because the projection $\Pi_0$ has finite rank, see also
the proof of Lemma~8 below).

\smallskip

If $|z| < 1/R$ is in the spectrum of $\MM_k$ on $\BB_{k}$,
the argument with mollifiers described above can be used 
to invoke  the ordinary Plemelj-Smithies formula  as in Lemma~4.4.2 of [Go],
to see that the determinant $\Det^\flat(\Id+\DD_{k,L}^{(r)}(z))$
has at most a  pole of order the algebraic multiplicity of
the eigenvalue at $z$, proving claim (2). \qed
\enddemo

\smallskip

We  next relate the
zeroes of the analytic continuation of the formal  determinant
$\Det^\flat(\Id+\DD_{k,L}(z))$, to the presence of an
eigenvalue $-1$ for the operators $\DD_{k,L}(z)$:

\proclaim{Lemma 7 (Zeroes of $\Det^\flat(\Id+ \DD_{k,L}(z))$)}
Assume Axioms 1, 2 and 3.
For $k=0, \ldots, n-1$ and all large enough $L$,
$\DD_{k,L}(z)$ extends holomorphically on $V_k$
to a family of operators on $\AA_{k+1,\BB_{k+1 }}$. The essential
spectral radius of each  $\DD_{k,L}(z)$ is  $< 1$.
If $\Det^\flat(\Id+ \DD_{k,L}(z))=0$
for $z \in V_k$, then $-1$ is an eigenvalue
of $\DD_{k,L}(z)$ on $\AA_{k+1,\BB_{k+1 }}$.
\endproclaim

\proclaim{Sublemma}
Let $A$ be a compact subset of $\real^n$.  Let $\widetilde \LL$
be a Hilbert-Schmidt operator on $L^2(A)$ written in kernel form
$
\widetilde \LL\varphi(x) =\int \widetilde \KK_{xy} \varphi(y) \, dy 
$.

Let $\LL$ be a bounded operator acting
on a Banach space $\BB$ of distributions over $A$
containing $C^\infty(A)$.
Assume that for small $\epsilon>0$
there are $C^\infty$ kernels
$\KK_{\epsilon,xy}$ and $\widetilde \KK_{\epsilon,xy}$
so that the Fredholm determinants of the
associated operators   coincide
$$
\Det (\Id-\lambda \LL_\epsilon)=
\Det (\Id-\lambda \widetilde \LL_\epsilon) \, ,
\forall \epsilon \, ,
$$
and assume also
$
\lim_{\epsilon\to 0}
\| \LL_\epsilon -  \LL \|_{\BB\to \BB}=0$, 
$\lim_{\epsilon \to 0}
|\widetilde \KK_{\epsilon,xy} -
\widetilde \KK_{xy} |_{L^2(K\times K)}=0 
$.

Then, for all $\lambda$ with
$1/|\lambda| < \rho_{ess}( \LL)$,
writing 
$$\Det_2^{\hbox{reg}}(\Id-\lambda \widetilde \LL)=\exp -\sum_{\ell =2}^\infty {\lambda^\ell\over \ell}
\tr \widetilde \LL^\ell \, , 
$$
we have
$
\Det_2^{\hbox{reg}}(\Id-\lambda \widetilde \LL)=0 \Longrightarrow \lambda^{-1}
\hbox{ is an eigenvalue of } \LL \hbox{ on }
\BB
$.
\endproclaim

\demo{Proof of the Sublemma}
The statements in this proof hold for all
$1/|\lambda| < \rho_{ess} (\LL)$ (uniformly in
any compact subset).
$\Det(\Id-\lambda \LL_\epsilon)=
\Det(\Id-\lambda \widetilde \LL_\epsilon)$ vanishes if and only if
$1/\lambda$ is an
eigenvalue of $\LL_\epsilon$ (on $L^2(A)$, or equivalently
on $C^\infty(A)$, using that the image of an element
of $L^2(A)$ by an operator with $C^\infty$ kernel is $C^\infty$), using the l.h.s.,
if and only if $1/\lambda$ is an
eigenvalue of $\widetilde  \LL_\epsilon$
(on $L^2(A)$ or equivalently on $C^\infty(A)$) using the r.h.s.
The convergence of the kernels 
$\widetilde \KK_{\epsilon,xy}$
implies both that $\Det_2^{\hbox{reg}}(\Id-\lambda \widetilde \LL_\epsilon)$
converges to
$\Det_2^{\hbox{reg}}(\Id-\lambda \widetilde \LL)$
(which is entire in $\lambda$), and
that $\widetilde \LL_\epsilon$ on $L^2(A)$ converges to
the compact operator $\widetilde \LL$ on $L^2(A)$.
Thus every zero $\lambda_0$ of
$\Det_2^{\hbox{reg}}(\Id-\lambda \widetilde \LL)$  is a limit of
$\lambda_\epsilon$ so that $1/\lambda_\epsilon$ is an
eigenvalue of
$\widetilde \LL_\epsilon$ (on $L^2(A)$).
(Indeed, such a zero corresponds to $1/\lambda_0$ being an
eigenvalue of $\widetilde \LL$.)
By the first observation, $1/\lambda_\epsilon$ is an
eigenvalue of  $\LL_\epsilon$
(on $C^\infty(A)$, so that the eigenfunction
is in $\BB$). If
$|\lambda_0| < 1/\rho_{ess} (\LL)$,
since $\LL_\epsilon$ converges to $\LL$
in operator
norm on $\BB$, if $1/\lambda_\epsilon$ is an
eigenvalue of $\LL_\epsilon$ for all $\epsilon$
then $\lambda_0$ is  an eigenvalue of $\LL$
(for an eigenvector which is the $\BB$-norm limit of
$\varphi_\epsilon\in C^\infty(A)$).
\qed
\enddemo

\demo{Proof of Lemma 7}
We have shown that
the left-hand-side of \thetag{3.13}
extends holomorphically to $V_k$. By definition and the Axioms
the operator $\DD_{k,L}(z)$  extends holomorphically
to $V_k$ in the sense of bounded operators on 
$\AA_{k+1, \BB_{k+1 }}$.  By the corollary
of Lemma ~5, it is  enough to show the claim about the zeroes
of the determinant.

Taking iterates to get Hilbert-Schmidt operators on $\AA_{k+1,L^2(K')}$,  one can
apply the sublemma to show that if $z \in V_k$
is such that $\Det^\flat(\Id+\DD_{k,L}(z))=
\Det^\flat(\Id+\DD^{(r)}_{k,L}(z))=0$
then $-1$ is an eigenvalue of $\DD_{k,L}(z)$ on 
$\AA_{k+1, \BB_{k+1 }}$. 

More precisely, for fixed $z\in V_k$, set
$
\LL = (\DD_{k,L}(z))^m  \, , \quad
\tilde \LL = (\DD_{k,L}^{(r)}(z))^m  
$.
Then, to construct the smooth kernels
required by the assumptions of the sublemma,
we (again) smoothen $\DD_{k,L}(z)$ by
pre and post-convolution with a $C^\infty$ mollifier 
$\delta_\epsilon$, 
writing $\DD^\epsilon_{k}(z)=
\delta_\epsilon \DD_{k,L}(z)\delta_\epsilon$ for the
new operator.
If we mollify $\DD^{(r)}_{k,L}(z)$  similarly, the equality
\thetag{3.13} between determinants remains true,
and the kernels of $\DD^{\epsilon,(r)}_{k}(z)$
converge as $\epsilon \to 0$, in
the $L^2(K\times K)$ topology. To see that
$\DD^\epsilon_{k}(z)$ converges to $\DD_{k,L}(z)$,
in the sense of operators on $\AA_{k+1, \BB_{k+1 }}$,
use  that $\delta_\epsilon$ converges to $\Id$  in the
$\AA_{k+1,\BB_{k+1 }}$ topology. 
\qed
\enddemo

\smallskip

{\bf A modification of the homotopy operators}

\smallskip

We next explain how to exploit Lemmas 6--7.
We first 
discuss the effect of making a suitable finite-rank
perturbation of $\SS$ in the definition
of the kneading operators:

\proclaim{Lemma 8 (Perturbing the homotopy operators)} Assume Axioms 1 and 2.
Let $\SS'$ be a finite-rank perturbation of $\SS$
(i.e., $\SS'_k-\SS_k$ is finite-rank from
$\AA_{k+1, \BB_{k+1}}$ to $\AA_{k,\BB_k}$),
so that $d\SS'+\SS' d =\Id$ and  $(\SS')^2 =0$. 
Then
the statements of Lemma ~1  and the modified 
equality  \thetag{3.10} from Theorem~2
remain true for  $\SS'$ and for the perturbed
kneading operators $\DD'_{k,L}(z)$ defined by
$$
\DD'_{k,L}(z)= \NN_{k,L}(z) (\Id- \MM_{k,L}(z))^{-1} \SS'_k \, , k=0,
\ldots, n-1 \, . \tag{3.16}
$$
\endproclaim

\demo{Proof of Lemma 8}
Since $d\SS'+\SS'd =\Id$, Baillif's  proof
extends relatively 
straightforwardly to  the modified convolution operators $\SS'$.
Indeed, we may use the decomposition
$$
\DD_{k,L}'(z)=\DD_{k,L}(z)+ ( \DD'_{k,L}(z)-\DD_{k,L}(z))
=\DD_{k,L}(z)+\FF_{k,L}(z)\, ,
$$
into a power series with coefficients operators whose
Schwartz kernel is well-behaved along the diagonal
(in particular the flat trace exists, and Fubini is allowed,
so that all desired commutations hold [Bai]), summed with a power series
whose coefficients are finite-rank operators.
These finite-rank operators act on $\AA_{k+1, \BB_{k+1 }}$
and each coefficient of $\DD_{k,L}(z)$ is bounded on this space.
Hence, we may define the formal traces
and determinants (as power series only) by setting
$$
\tr^\flat (\DD_{k,L}(z)+\FF_{k,L}(z))^m =
\tr^\flat \DD_{k,L}(z)^m + \tr^\flat \FF_{m,k,L}(z)\, ,
$$
where $\FF_{m,k,L}(z)$ is a power
series with coefficients finite-rank operators on $\AA_{k+1, \BB_{k+1 }}$,
for which the flat trace is just defined to be the sum of
eigenvalues.

In fact, since $\DD'_{k,L}(z)-\DD_{k,L}(z)$ starts with a precomposition
by $\chi_K(\SS'_k-\SS_k)$, which is finite-rank from
$\AA_{k+1, \BB_{k+1 }}$ to $\AA_{k, \BB_{k }}$,
the resolvent factor $(\Id-\MM_{k,L}(z))^{-1}$ is well defined
if $z\in V_k$. Thus, the power series $\tr^\flat \FF_{m,k,L}(z)$
defines in fact  a holomorphic function in $V_k$. (This argument
will be useful in the proof of Lemma~10 below.)
\qed
\enddemo


\proclaim{Lemma 9 (Adapted homotopy operators $\SS'$)}
Assume  Axioms ~1 and~  2, and
suppose that all eigenvalues of modulus $> R$ of the $\MM_k$
on $\AA_{k,\BB_k}$($k \le n-1$)
are simple, and that the corresponding eigenvectors $\varphi_k$ 
do {\it not} belong to
$\Ker d=\Imm d$ or to $\Ker \NN_k$.
Then there are finite-rank operators
$$
\eqalign{
&\FF_{k}: \AA_{k+1, \BB_{k+1}}\to
\{ \varphi \in \AA_{k} \mid \chi_K \varphi \in \AA_{k, \BB_{k }} \} \, ,
\quad 0 \le k \le n-1 \, ,\cr
&
\hbox{so that the perturbed operators }\qquad
\SS'_k = \SS_k +\FF_k\, ,}
$$  
satisfy $\SS' \SS'=0$, $\SS'd+ d\SS'=1$ and, additionally,
if $0 \le k \le \min(s,n-1)$ and $1/z$ with $|z|> R$ is an eigenvalue  of
$\MM_k$ on $\AA_{k, \BB_{k }}$ with corresponding (simple)
fixed vector $\varphi_k(z)$ for $\MM_{k,\ell}(z)$ and dual
fixed vector $\nu_k(z)$, with $\nu_k(z)(\varphi_k(z))=1$, then
for large enough $\ell$ 
$$
|\nu_k (z) (\chi_K \SS'_k  \NN_{k,\ell}\varphi_k(z))|  > 0 \, .
$$
\endproclaim

\demo{Proof of Lemma 9}
For each $k$ we consider the finite set of
``bad pairs''  $(z, \varphi_k)$, with
$z=z_{k,i}\in \complex$ and $\varphi_{k,i}$ an fixed vector such that 
$\nu_{k,i} (z) (\chi_K \SS_k \NN_{k,\ell} \varphi_{k,i}(z))= 0$.
We write the argument assuming that this set is either empty or
a singleton 
$\{(z, \varphi_k)\}$ for each
$k$ in order to simplify notation. Our assumptions imply
that $\upsilon_k=\nu_k(\SS_k d \varphi_k)\ne 0$ and 
$\NN_{k,\ell} (z) \varphi_k(z)\ne 0$ for all $\ell$.
For each $k$ so that a bad pair $(z,\varphi_k)$ exists,
set $\varphi'_{k}=\varphi_k$ and
let $\nu'_{k}$, $\nu'_{k+1}$  be   continuous
functionals on $\AA_{k , \BB_{k }}$,
$\AA_{k+1, \BB_{k+1}}$, respectively, which satisfy 
$\nu'_k \circ d = 0$, $\nu'_{k+1} \circ d = 0$, and 
$$
\alpha_{k,\ell}=\nu'_{k+1}(\NN_{k,\ell} \varphi_k) =-\nu'_{k+1}(\MM_{k+1,\ell}d \varphi_k)\ne 0
\, , \forall \ell \, , \quad \beta_k=\nu'_k(\varphi_k)\ne 0\, .
$$
If there is no $\varphi_k$ we set $\varphi'_{k}=0$ and $\nu'_{k+1}=0$
and allow $\alpha_k$ and $\beta_{k,\ell}$ to vanish.
We put, for complex $\epsilon$ of small modulus, and all $k$
$$
\FF'_{k} (\varphi)= -\epsilon d \varphi'_{k-1} \cdot \nu'_{k}(\SS_k \varphi)
+\epsilon \varphi'_{k} \cdot \nu'_{k+1}(\varphi)\, .
$$
We have $d\FF'_{k-1}=-\FF'_k d$ so that $d(\SS+\FF')+(\SS+\FF')d=\Id$. 
Thus, setting $\FF=\FF'-d\SS\FF'-d\FF'\SS-d\FF'\FF'$ we
have $d\SS'+\SS'd=\Id$ and $\SS'\SS'=0$. 
Finally, for each $k$ with a bad pair, we find for
uncountably many small values of $\epsilon$ and
all large enough $\ell$
$$
\eqalign
{
\nu_k((\SS'_{k}\NN_{k,\ell})  \varphi_k)&=
\epsilon  
\biggl (\nu_k(\SS_k d \varphi_{k}) \cdot \nu'_{k+1} (\NN_{k,\ell}(z)\varphi_k)  \cr
&\quad-
\nu_k(d \varphi_{k -1})\cdot 
\bigl [\nu'_{k}(\SS_k \NN_{k,\ell}(z)\varphi_k)
+\epsilon \nu'_k(\varphi_k)\nu'_{k+1} (\NN_{k,\ell}(z)\varphi_k) 
\bigr ] \biggr )\cr
&=  \epsilon \bigl (\upsilon_k \alpha_{k,\ell}+ \nu_k(d\varphi_{k-1})
(\nu'_{k}(\SS_k \NN_{k,\ell}(z)\varphi_k)+\epsilon \beta_k \alpha_{k,\ell}) 
\bigr )\ne 0  \, . 
}
$$
By taking $\epsilon$ small enough we may ensure that no new bad pairs are
created. \qed
\enddemo

\noindent Replacing $\SS$ by the finite rank perturbation
$\SS'$,  we get a stronger
version of Lemmas~6--7:

\proclaim{Lemma 10 (Meromorphic extension -- guaranteeing poles)}
Under the assumptions of Lemma ~9, and up to taking
a larger value of $L$, for each $k=0, \ldots, n-1$, Axiom~3 implies:

\noindent (1)
The power series $\Det^\flat(\Id+ \DD'_{k,L}(z))$ 
defines a holomorphic function in $V_k$.

\noindent (2)
If  $|z_0| < 1/R$  is such that $1/z_0$ is a simple eigenvalue
of $\MM_k|_{\AA_{k, \BB_{k }}}$,
then $\Det^\flat(\Id+\DD'_k(z))$ is meromorphic
at $z_0$, with a pole of order exactly one.

\noindent (3)
$\DD'_{k,L}(z)$ extends holomorphically on $V_k$
to a family
of operators on $\AA_{k+1, \BB_{k+1}}$,
each such $\DD'_{k,L}(z)$ has essential
spectral radius strictly less than $1$,
and $\Det^\flat(\Id+ \DD'_{k,L}(z))=0$ 
if and only if $\DD_{k,L}'(z)$ has an eigenvalue $-1$ on 
$\AA_{k+1,\BB_{k+1 }}$.
\endproclaim

\demo{Proof of Lemma 10}
From the
proof of Lemma~6, $\DD^{(r)}_k(z)$ has its $[n/2]+1$th iterate
trace-class on $\AA_{k+1, L^2(K')}$ for 
all $z \in V_k$.
Using Axioms~1--2, for each $z\in V_k$, 
$$
\eqalign
{\DD^{(r)'}_{k,L}(z)
&=\DD^{(r)}_{k,L}(z)+
(\widetilde \JJ_{k})^{-1} \NN_{k,L}(z)(\Id-\MM_{k,L}(z))^{-1}
(\SS'_k-\SS_k)\widetilde \JJ_{k} 
}
$$ 
is such that its $[n/2]+1$th iterate is a finite-rank perturbation on
$\AA_{k+1, L^2(K')}$ of $\DD^{(r)}_{k,L}(z)$.
The sum of the $L^2(K')$ (or, equivalently, flat)
trace of $(\DD^{(r)}_{k,L}(z))^{[n/2]+1}$ and the flat trace of
a composition (in any order) of $j$ factors of the
finite-rank term $[\DD^{(r)'}_{k,L}(z)-\DD^{(r)}_{k,L}(z)]^{[n/2]+1}$  with
$([n/2]+1-j)$ factors $\DD^{(r)}_{k,L}(z)$ is holomorphic in $V_k$.
Its power series is equal
to the sum of the formal flat trace of $\DD_{k,L}^{[n/2]+1}(z)$
with the trace of a finite rank operator
which has the same trace as $(\DD'_{k,L}(z)-\DD_{k,L}(z))^{[n/2]+1}$.
(Recall 
that $\DD'_{k,L}(z)-\DD_{k,L}(z)$ starts with a precomposition
by $\chi_K(\SS'_k-\SS_k)$ which is finite-rank from 
$\AA_{k+1, \BB_{k+1}}$ to $\AA_{k, \BB_{k}}$.)
Therefore $\tr^\flat  (\DD^{(r)'}_{k,L}(z))^{[n/2]+1}=\tr^\flat (\DD'_{k,L}(z))^{[n/2]+1}$,
which yields the first claim of Lemma~10.

Next, we may apply the Sublemma, essentially as in Lemma~7, to see that
the zeroes of the determinant $\Det^\flat(\Id+ \DD'_{k,L}(z))$ correspond to
eigenvalues $-1$, using  
$$
\eqalign{
\DD^{\epsilon, (r')}_{k,L}(z)&=
\DD^{\epsilon, (r)}_{k,L}(z)
+ \delta_\epsilon (\widetilde \JJ_{k})^{-1}\NN_{k,L}(z)(\Id-\MM_{k,L}(z))^{-1}
(\SS'_k-\SS_k)\widetilde \JJ_{k}\delta_\epsilon \cr
\DD^{\epsilon'}_{k,L}(z)&=
\DD^\epsilon_{k,L}(z)
+ \delta_\epsilon \NN_{k,L}(z)(\Id-\MM_{k,L}(z))^{-1}
(\SS'_k-\SS_k)\delta_\epsilon\, , }
$$
which both have a  finite-rank
(in $L^2(\real^n)$) second term.   This gives  Lemma~10,\therosteritem{3}.

If $1/z_0$ is a simple eigenvalue for $\MM_k$ with $z_0\in V_k$
then 
one proves, like in Lemma~ 6, that $\Det^\flat(\Id+\DD'_{k,L}(z))$
is meromorphic at $z_0$.
 with a pole of order at most one. 
Finally, 
we shall prove that $\Det^\flat(\Id+\DD'_{k,L}(z))$
does not have a removable singularity at $z_0$
showing that the pole has order exactly one 
i.e., claim
\therosteritem{2}.

For this, we use a
spectral decomposition of $\MM_{k,\ell}(z)$ for the
simple eigenvalue $\lambda_z$ on $\AA_{k,\BB_{k}}$,
for $z$ close to $z_0$ (with $\lambda_{z_0}=1$):
$$
(\Id-\MM_{k,\ell}(z))^{-1} = {\varphi_{k}(z) \over \lambda_z-1} \cdot
\nu_k(z) +
\RR_k(z)\, , \tag{3.17}
$$
with $\lambda_z$ holomorphic in $z$;
$\RR_k(z): \AA_{k,\BB_{k}} \to \AA_{k,\BB_{k}}$ depending
holomorphically on $z$ at $z_0$;
$\RR_k (z) \varphi_k(z) = \nu_k(z) \RR_k(z) = 0$;
the eigenvector $\varphi_k(z)  \in \AA_{k,\BB_{k}}$
of unit norm  depending holomorphically on
$z$ with $\varphi_k(z_0) =\varphi_k \in \AA_{k,\BB_{k}}$ so that $\MM_k(z_0) \varphi_k =\varphi_k$;
and  $\nu_k(z)$ unit-norm linear functionals  on $\AA_{k,\BB_{k}}$,
depending holomorphically on $z$,
with $\nu_k(z)(\varphi_k(z)) = 1$.

By Lemma 9, 
$|\nu_k(z_0) (\chi_K \SS'_{k} \NN_{k,\ell}(z_0) \varphi_k(z_0)) | > 0$
if $\ell$ is large enough.
It is then easy to see that
$
(z-z_0)\Det^\flat(\Id+\DD'_{k,L}(z))
$
does not vanish at $z_0$. 
\qed
\enddemo

\smallskip
{\bf Using the new homotopy operators}
\smallskip

Making use of the  homotopy operators from Lemma~9, we  
state the final ingredients needed in
our main result. The eigenvalues of modulus larger than $R$ of the
operators $d \SS \MM_k$  produce zeroes of the flat
kneading determinants (Lemma~11). To show that
such eigenvalues  do not contribute to
the zeroes and poles of the sharp determinant (except when
$\MM_k \tilde\varphi_k=d\SS\MM_k\tilde \varphi_k$,
in particular if $k=n$), we shall prove that they are not
intrinsic. More precisely, in Lemma~12  we 
construct perturbed homotopy operators
$\SS''$ which cause these eigenvalues to vary.

\proclaim{Lemma 11 (Zeroes of the flat kneading determinants)}
Under the assumptions of Lemma~9, and up to taking a larger
value of $L$,
for all $k \in \{1, \ldots, n\}$ and $|z| < 1/R$, 

\noindent (1)
If
$\Det^\flat(\Id+ \DD'_{k-1,L}(z))= 0$  then
$1/z$ is an eigenvalue of $d\SS'_{k-1}\MM_k$ acting on $\AA_{k,\BB_k}$
or $k\ge 2$ and $1/z$ is an eigenvalue of $d\SS'_{k-2}\MM_{k-1}$ 
acting on $\AA_{k-1,\BB_{k-1}}$, in either case, the geometric
multiplicity of the eigenvalue is at least equal to
the order of the zero.

\noindent (2)
If $1/z$ (with $|z| < 1/R$)
is an eigenvalue of $\MM_n$ 
acting on $\AA_{n,\BB_n}$ and is not an eigenvalue
of $\MM_{n-1}$ 
acting on $\AA_{n-1,\BB_{n-1}}$,
then $\Det^\flat(\Id+ \DD'_{n-1,L}(z))= 0$ and the order
of the zero is at least the algebraic multiplicity of
the eigenvalue.

\noindent (3)
If $k \ge 2$ and $1/z$ (with $|z| < 1/R$)
is an eigenvalue of $d\SS'_{k-2}\MM_{k-1}$ 
acting on $\AA_{k-1,\BB_{k-1}}$ 
then $\Det^\flat(\Id+ \DD'_{k-1,L}(z))= 0$ and the order
of the zero is at least the geometric multiplicity
of the eigenvalue (which is thus finite and the eigenvalue isolated).
\endproclaim

\demo{Proof of Lemma 11}
By Lemma~10, the assumption $\Det^\flat(\Id+ \DD'_{k-1,L}(z))= 0$
implies that   $1/z \notin \sp(\MM_{k-1})$ and
there is a nonzero $\varphi\in  \AA_{k,\BB_{k}}$  with
$\DD'_{k-1,L}(z)\varphi=-\varphi$. It follows that $\varphi$ is {\it not}
in the kernel of $\SS'_{k-1}$.
Writing $\varphi=d\SS'_{k-1}\varphi + \SS'_k d\varphi=\varphi_1+\varphi_2$ (with
$\SS'_{k-1}\varphi_1\ne 0$), our assumption implies
$$
\cases
d\SS'_{k-1}\varphi_1 +  d \SS'_{k-1}\NN_{k,L}(z)(\Id-\MM_{k-1,L}(z))^{-1}
\SS'_{k-1} \varphi_1=0&\cr
\varphi_2=\SS'_k d\varphi_2 =
-\SS'_k d \NN_{k,L}(z) (\Id-\MM_{k-1,L}(z))^{-1} \SS'_{k-1} \varphi_1 \, .&\cr
\endcases \tag{3.18}
$$
The first equality in \thetag{3.18} is equivalent with
$$
\eqalign
{
0&=\left ( d+ (d\MM_{k-1,L}(z) -d\SS'_{k-1} \MM_{k,L}(z) d)
(\Id-\MM_{k-1,L}(z))^{-1} \right )\SS'_{k-1} \varphi_1\cr
&=(\Id-d\SS'_{k-1} \MM_{k,L}(z))
d (\Id-\MM_{k-1,L}(z))^{-1} \SS'_{k-1} \varphi_1 \, .}
$$
Note that  $d (\Id-\MM_{k-1,L}(z))^{-1} \SS'_{k-1} \varphi_1$
is in $\AA_{k, \BB_k}$ by the boundedness of the
$d$ operator in Axiom ~1. 
Then it is not very difficult to check (see the proof
of Lemma~5) that the bounds on the norm
of $\NN_{k-1}^{(j)}$  translate into exponentially decaying
bounds for $j\mapsto \|\NN_{k-1,L}(z)^{(j)}\|$.

Since $\chi_K\SS'_{k-1}\varphi_1\in  \AA_{k-1, \BB_{k-1}}$
does not vanish,
$(\Id-\MM_{k-1,L}(z))^{-1} \SS'_{k-1} \varphi_1\ne 0$. If
$$
\hat \phi=d(\Id-\MM_{k-1,L}(z))^{-1} \SS'_{k-1} \varphi_1\ne 0 \, ,\tag{3.19}
$$
then we are done.
Indeed,
$\hat \phi\in \Imm d\cap \AA_{k,\BB_k}$ would then be a nonzero fixed vector
for   $d\SS'_{k-1} \MM_{k,L}(z)$. Since the essential
spectral radius of  $d\SS'_{k-1} \MM_{k,L}(z)$ is smaller
than one, is an easy algebraic exercise
to see that the fixed vectors $d\SS'_{k-1} \MM_{k,L}(z)$
(if $|z|< 1/R$) are in bijection with the eigenvectors of
$d\SS'_{k-1} \MM_{k}$ for the eigenvalue $1/z$ so that
we are in the first case of the first claim.
Note also at this point that if there exists $\hat\phi \in \AA_{n,\BB_n}$
with $d\SS'_{n-1}\MM_{n,L}(z) \hat\phi=\hat\phi$ then 
$\varphi_1:= d(\Id-\MM_{n-1,L}(z)) \SS'_{n-1}
\hat\phi\in \AA_{n,\BB_n}$ would satisfy the first identity in
\thetag{3.18}, while $\varphi_2$ may be defined by the second equality
of \thetag{3.18}. In this case $\Det^\flat(\Id+ \DD'_{n-1,L}(z))= 0$
by Lemma~10.

Let us assume that \thetag{3.19} does not hold. Writing
$\upsilon=(\Id-\MM_{k-1,L}(z))^{-1} \SS'_{k-1}\varphi_1$
we would then have $ \upsilon \in  \AA_{k-1, \BB_{k-1}}$  
and  $d\upsilon=0$ and,
applying $\SS'_{k-1} d$ to both sides of
$\SS'_{k-1}\varphi_1=d\SS'_{k-2}\upsilon-\MM_{k-1,L}(z) \upsilon$,
$$
\SS' _{k-1}\varphi_1=-\SS'_{k-1} d \MM_{k-1,L} (z)\upsilon\, ,
$$
since $\SS'_{k-1}\varphi_1=(\Id-\MM_{k-1,L}(z)) \upsilon$, we
find    $(\Id-\MM_{k-1,L}(z))\upsilon=-\SS'_{k-1}d \MM_{k-1,L}(z)\upsilon$,
contradicting  $1/z \notin \sp \MM_{k-1}$ if $k=1$, and 
if $k\ge 2$ implying
$$
d\SS'_{k-2}\upsilon=\upsilon= (\Id-\SS'_{k-1} d)\MM_{k-1,L} (z)\upsilon=
 d \SS'_{k-2} \MM_{k-1,L} (z)\upsilon \, . \, 
$$
For the converse (i.e., the last claim in
Lemma~11), we see that if $\upsilon= d \SS'_{k-2} \MM_{k-1,L} (z)\upsilon$,
we may take $\varphi_1=-d\MM_{k-1,L}\upsilon$ since then 
$(\Id-\MM_{k-1,L})^{-1}\SS'_{k-1}\varphi_1=-\upsilon$ so that
$d(\Id-\MM_{k-1,L})^{-1}\SS'_{k-1}\varphi_1=0$.
\qed
\enddemo

\smallskip

\proclaim{Lemma 12 (Modified homotopy operators $\SS''$: perturbing nonintrinsic eigenvalues)}
Assume that for some
$ 2 \le k \le n-1$ there are $0\ne \tilde \varphi_k \in \AA_{k, \BB_{k }}$ and
$z \in \complex$ with  $|1/z| > R$ such that
$$
\tilde \varphi_k = {1 \over z}d \SS'_{k-1} \MM_k\tilde \varphi_k \, ,
$$
as an isolated eigenvalue of finite geometric multiplicity.
Then, there are $z'\ne z$, arbitrarily close to $z$, and two
rank-one operators of arbitrarily small norm
$$
\FF'_{\ell}: \AA_{\ell+1, \BB_{\ell+1}}\to
\{ \varphi \in \AA_{\ell} \mid \chi_K \varphi \in \AA_{\ell, \BB_{\ell }} \} \, ,
\, \ell=k-1, k\, ,
$$
so that the perturbed operators
$
\SS''_{k-1}=\SS'_{k-1}+\FF'_{k-1}$, 
$\SS''_{k}=\SS'_{k}+\FF'_{k}$,
$\SS''_\ell = \SS'_\ell$,
$\ell\notin \{ k-1 , k\}
$,
still
satisfy $\SS'' \SS''=0$, $\SS''d + d\SS''=1$ and, additionally,
$1/z'$ is an eigenvalue of
$d \SS''_{k-1} \MM_k\varphi$ on $\AA_{k, \BB_{k }}$
while $d \SS''_{\ell} =d\SS'_\ell$ for all $\ell\ne k-1$.
\endproclaim

\demo{Proof of Lemma 12}
We take $\varphi'_{k-1} := \SS'_{k-1}\MM_k \tilde \varphi_k\in 
\AA_{k-1, \BB_{k-1}}$. By our assumptions,
$d\varphi'_{k-1}\ne 0$. Let then $\nu'_{k}$ be a unit-norm continuous
functional on $\AA_{k, \BB_{k}}$ which satisfies
$$
\nu'_k \circ d = 0 \, , \quad \alpha:=\nu'_k(\MM_k \tilde \varphi_k)\ne 0 \, .
$$
(Take  $\nu'_k (\varphi)= \nu'_{k+1}(d\varphi)$ 
with
$\nu'_{k+1}$ continuous  on $\AA_{k+1, \BB_{k+1}}$ and
$\nu'_{k+1}(d\MM_k \tilde \varphi_k)\ne 0$.)
We set, for small complex $\epsilon$,
$$
\eqalign
{
\FF'_{k-1} (\varphi)= \epsilon \varphi'_{k-1} \cdot \nu'_{k}(\varphi) \, , \quad
\FF'_{k} (\varphi)= -\epsilon d \varphi'_{k} \cdot \nu'_{k}(\SS_k \varphi)\, .
}
$$
We have $d\FF'_{k-1}=-\FF'_k d$, $\FF'_{k-1}d=0$ and
$d \FF'_k=0$ so that $d\SS''+\SS''d=\Id$. Also,
$\FF'_k \SS'_{k+1}=0$, and 
$\SS'_{k-1}\FF'_k + \FF'_{k-1} \SS'_k + \FF'_{k-1}\FF'_k =0$, guaranteeing
$\SS''\SS''=0$. Clearly, $d\FF'_k =0$ so that $d\SS''_k=d\SS'_k$.
Finally, 
$$
z d \SS''_{k-1} \MM_k \tilde \varphi_k=
\tilde \varphi_k + \epsilon z \cdot (d\SS'_{k-1}\MM_k \tilde \varphi_k)
\cdot \nu'_k (\MM_k \tilde \varphi_k) =  (1+ \alpha \epsilon) \tilde \varphi_k \, . 
\hbox{\qed}
$$
\enddemo


\pagebreak

{\bf Main Result}

\smallskip

\proclaim{Theorem 13}
Let $\psi_\omega$, $g_\omega$ satisfy the
assumptions of Section ~2,  let $\BB_{k,t}$ and $R$ satisfy
Axioms~ 1 and ~2--3.
Then
$\Det^\#(\Id-z\MM)$ is meromorphic in the disc $\{ |z| < 1/ R\}$.
The order
of $z$ as a zero/pole of $\Det^\#(\Id-z\MM)$
coincides with the sum of the algebraic multiplicity
of $1/z$ as an eigenvalue of the  $\MM_{2k}$ on $\AA_{2k, \BB_{2k}}$,
for $0 \le 2k \le n$, minus the sum of the algebraic multiplicity
of $1/z$ as an eigenvalue of the $\MM_{2k+1}$ on $\AA_{2k+1, \BB_{2k+1}}$ for
$1 \le 2k+1 \le n$.
\endproclaim

\demo{Proof of Theorem 13}
Perturbing  our family,  we may assume that all eigenvalues
of moduli $> R$ of the $\MM_k$,
$k=0, \ldots, n-1$, are simple,  and none of their eigenvectors
belong to $\Ker d$ or $\Ker \NN_k$. We may also assume that
the eigenvalues of $\MM_n$ are simple.

We apply Lemma ~9 to construct adapted homotopy operators $\SS'$
and let $\DD'_{k,\ell}(z)$  be the kneading
operators from \thetag{3.16}, for large enough $\ell\ge L$.

We consider the modified equality \thetag{3.10} from Theorem~2
for the perturbed family
and the adapted homotopy operators,
using  Lemma~10 to view it  as an
alternated product of meromorphic
functions in the disc of radius $1/R$.
(Later in the proof we shall use further finite rank perturbations of  the $\SS'$
given by Lemma ~12 and thus
satisfying the assumptions of Lemma~8.)

By Lemma 10,  the determinant
$\Det^\flat (\Id+\DD'_{2k,\ell}(z))$  for $0\le 2k\le n-1$
(which appears in
the denominator) has poles in the disc of radius $1/R$ only at the inverse
eigenvalues of $\MM_{2k}$, with order exactly one.
We have the same statement for the
$\Det^\flat (\Id+\DD'_{2k+1,\ell}(z))$  for $1\le 2k+1\le n-1$,
which appear in the numerator.

Lemma 11 also says that if   $\Det^\flat (\Id+\DD'_{k-1,\ell}(z))$ for
$1\le k \le  n$ (in the denominator for even $k-1$
and in the numerator for odd $k-1$) vanishes in this disc 
then $1/z$ is an eigenvalue
of $d \SS'_{k-1}\MM_{k}$ or $k \ge 2$ and
$1/z$ is an eigenvalue of $d \SS'_{k-2}\MM_{k-1}$.
Also, whenever
$1/z$ is an eigenvalue of $d \SS'_{k-2}\MM_{k-1}$ 
then $\Det^\flat (\Id+\DD'_{k-1,\ell}(z))=0$, and if $1/z$ is an eigenvalue
of $d\SS \MM_n=\MM_n$ then $\Det^\flat (\Id+\DD'_{n-1,\ell}(z))=0$.

To finish the proof, we will show that the zeroes of the
flat determinants must cancel in the alternated product,
except of course for the  zeroes of $\Det^\flat (\Id+\DD'_{n-1,\ell}(z))$ 
(in the denominator if $n-1$ is even,
in the numerator if $n-1$ is odd) corresponding to 
$1/z$ being an eigenvalue of $\MM_n$.
Assume for a contradiction that  $\Det^\#(1-\hat z\MM)= 0$ to order
$D$ strictly larger than the value claimed in Theorem~13,
due (at least in part) to a factor $\Det^\flat (\Id+\DD'_{\hat k-1,\ell}(\hat z))=0$
vanishing to order $\hat D\ge 1$
for some odd  $0\le \hat k-1\le n-1$. (The case of poles and even $k-1$
is dealt similarly.)
Lemma~11 says that $1/\hat z$ is then  either an eigenvalue
of $d\SS'_{k-2}\MM_{\hat k-1}$ (since $\hat k\ge 2$) or an eigenvalue of
$d\SS'_{\hat k-1}\MM_{\hat k}$ (and $\hat k < n$). If  $1/\hat z$ is an eigenvalue
of $d\SS'_{k-2}\MM_{\hat k-1}$ (the other case is left to the reader), 
we may use Lemma~12 
to perturb this eigenvalue to some
$1/\hat z'\ne 1/\hat z$. Note that this does
not modify  the $\Det^\flat (\Id+\DD'_{k-1,\ell}( z))$ 
except for $k= \hat k$ and $k=\hat k+1$. (For both of  these determinants, 
only the set of zeroes may change, and only one of the determinants is in
the numerator.) 
The perturbation $\DD''_{\hat k-1,\ell}(z)$ of $\DD'_{\hat k-1,\ell}(z)$
may be made as small as desired in (finite-rank) operator norm in  a
neighbourhood of $\hat z$
by taking small enough $\epsilon\ne 0$ in Lemma ~12. Since
$\hat z\ne \hat z'$, the continuous dependence
of the regularised determinant on the operator together with Rouch\'e's Theorem guarantee
that the order of $\hat z$ as a zero of $\Det^\flat (\Id+\DD''_{\hat k-1,\ell}(z))$
is strictly smaller than $\hat D$. Iterating this procedure (or the procedure
associated to the other
case) at most $\hat D$ times, we get that $\Det^\flat (\Id+\DD'''_{\hat k-1,\ell}(\hat z))\ne 0$,
while none of  the zeroes of the other kneading determinants (in the numerator)
or the poles (in the denominator or in fact also the numerator) have been altered.
Then, either the order of $\hat z$ as a zero
of $\Det^\#(1- z\MM)$ is strictly smaller than $D$, a contradiction,
or $\Det^\flat (\Id+\DD'''_{k-1,\ell}(\hat z))=0$ for some odd $k-1\ne
\hat k-1$. In the second case, we may proceed as above to ensure
$\Det^\flat (\Id+\DD''''_{k-1,\ell}(\hat z))\ne 0$ and eventually obtain
the contradiction that  the order of $\hat z$ as a zero
of $\Det^\#(1- z\MM)$ is strictly smaller than $D$.
\qed
\enddemo

\noindent A corollary of the proof of Theorem~13 is
(we do not have to replace
$\widetilde R$ by $\widetilde R^{1/n}$):

\proclaim{Theorem 14}
Assume Axioms~1 and 2. For each large enough
$L$, the following alternated product
of regularised determinants  extends holomorphically
to $\{ |z|< 1/\widetilde R\mid 1/z \notin \cup_k \sp(\MM_k)\}$:
$$
\prod_{k=0}^{n-1} \Det^{\hbox{reg}}_{[n/2]+1}(\Id+\DD_{k,L}^{(r)}(z))
^{(-1)^{k+1}} \, .
\tag{3.20}
$$
\endproclaim


Although
$
\Det^{\hbox{reg}}_{[n/2]+1}(\Id+\DD_{k,L}^{(r)}(z))
=\exp -\sum_{\ell = [n/2]+1}^\infty {z^\ell\over \ell}\tr^\flat (\DD_{k,L}(z))^\ell
$,
it does not seem easy to relate \thetag{3.20} to a dynamical
zeta function.

\smallskip

\head
4. Application to expanding maps on compact manifolds
\endhead

In this section, we  apply Theorem~13 from Section~ 3 to
$C^{r}$ (locally) expanding endomorphisms on compact manifolds
and $C^r$ weights,
giving a new (and completely different) proof of a
result of Ruelle [Ru3], avoiding Markov structures.  We do not  recover
the full strength of his statement except if the dynamics
is $C^\infty$ and if the weight is the inverse Jacobian
(with respect to Lebesgue).

\smallskip

Let us start by stating precisely this result.
Let $M$ be a $C^{r}$ ($r\ge 1$) compact manifold of dimension $n\ge 2$, and let
$f:M\to M$ be $C^{r}$ and (locally) uniformly  expanding,
that is, there is $\theta <1$ such that
$\|D_a f(\xi)\|\ge\theta^{-1}\|\xi\|$  for all $\xi$ in the tangent space $T_a M$
at $a$, where $\|\cdot\|$ is the Euclidean norm on $T_a M$.
Let $g:M\to\complex$ be  $C^r$.
The transfer operator $\LL_0=\LL_{f,g}$, acting on
the Banach space of $C^r$ functions $M\to\complex$,
is given by
$$
\LL_0\phi(a)=\sum_{b:f(b)=a}g(b)\phi(b) 
=\sum_j g_j(a) \phi(\psi_j (a)) \, ,
\tag{4.1}
$$
where the $\psi_j$ are the finitely many local inverse branches of $f$.
Similarly, we can introduce operators acting on Banach spaces of
$k$-forms on $M$ with $C^m$ coefficients, for $0\le m \le r-1$, putting
$
\LL_k\phi=\sum_{j}g_j \cdot (\psi_j^*) \phi 
$.

An expanding map is transversal, and the Lefschetz numbers of
the inverse map at the periodic orbits are all positive,
so that the Ruelle zeta function associated to $f$,~ $g$
$$
\zeta_{f,g}(z)=\exp \sum_{m=1} ^\infty{z^m\over m}
\sum_{a \in \Fix f^m} \prod_{\ell=0} ^{m-1}  g(f^\ell(a))
$$
can be viewed as a Lefschetz-Ruelle zeta function.
Ruelle proved:

\proclaim{Theorem (Ruelle [Ru3])}
Let $P\in \real$ be the topological pressure of $\log|g|$ and $f$.

\noindent (1)
The spectral radius of $\LL_0$ acting on
$C^m$ is at most $e^P$ while its the essential
spectral radius is at most $\theta^{m}e^P$ for  $0\le m \le r$.
The spectral radius of $\LL_k$ acting on
forms with $C^m$ coefficients is at most $\theta^{k}e^P$ while its the essential
spectral radius is at most $\theta^{m+k}e^P$, for 
$1 \le k\le n$ and $0\le m \le r-1$.

\noindent (2)
The power series  $\zeta_{f,g}(z)$  defines a meromorphic function
in the disc of radius $\theta^{-r}e^{-P}$.
In this disc, the order
of $z$ as a zero/pole of $\zeta_{f,g}(z)$
coincides with the sum of the algebraic multiplicity
of $1/z$ as an eigenvalue of the  $\LL_{2k+1}$ acting on
$C^m$ for any $r-(2k+1)\le m \le r-1$ and $1 \le 2k+1 \le n$, minus the sum 
of the algebraic multiplicity
of $1/z$ as an eigenvalue of the $\LL_{2k}$ acting on $C^m$  for $r-2k\le m \le r-1$ and
$0 \le 2k \le n$ (for $k=0$ one can also let $\LL_0$ act
on $C^r$).
\endproclaim

{\it We will recover Ruelle's result  if the system is in fact
$nr$ times differentiable, with the same estimates
if $g(y)=1/| \det Df (y)|$,
(i.e.,  the weight giving rise to the absolutely continuous invariant measure),
and only a weaker result, replacing $e^P$ by
$\theta^{-n+1} e^P$, in the case of an arbitrary smooth weight.} There
is certainly room for improvement here. 

\smallskip
Once the first claim of the above theorem is established
(see [GuLa] for better estimates),
Ruelle  [Ru3] associates to each $\LL_k$ a Fredholm-like
(flat) dynamical determinant
$$
d_k(z)=\Det^\flat (\Id - z\LL_k)= \exp -\sum_{m=1} ^\infty{z^m\over m}
\sum_{a \in \Fix f^m} \prod_{\ell=0} ^{m-1}  { g(f^\ell(a))\cdot  \hbox{Tr}\,  \Lambda^k (D_a f^{-m} )
\over \det(\Id - D_a f^{-m}) }\, ,\tag{4.2}
$$
and proves  that $d_k(z)$ is holomorphic in the disc of radius $\theta^{-k-r}e^{-P}$,
where its zeroes
correspond exactly the the inverse eigenvalues of $\LL_k$
(outside of the disc of radius $\theta^{k+r}e^{P}$). 
Writing   $\zeta_{f,g}(z)$ as an alternated product of the
$d_k(z)$ gives the second claim.

The present approach does not allow us to analyze the independent
factors $d_k(z)$.
However, since the spectral radii of the operators $\LL_k$ are strictly
decreasing,  the annulus
$e^{-P}\le|z|<\theta^{-1}e^{-P}$, e.g.,
only contains inverse eigenvalues of $\LL_0$ acting on $C^r(M)$.

\smallskip
{\bf From a manifold to $\real^n$ -- Equivalent models}
\smallskip

We wish to associate to $f$ and $g$
data $\{\psi_\omega,g_\omega\}_{\omega\in\Omega}$ in such a way that the operators $\LL_k$
and $\MM_k$ are conjugated. Their spectra on suitable spaces will thus coincide.
For this, first  choose a $C^r$ atlas $\{V_j\}$ for $M$
such that $f|_{V_j}\to f(V_j)$ is a diffeomorphism.
Choose a $C^r$ partition of the unity
$\{\chi_j\}$, where each $\chi_j$ is supported in $V_j$.
Denote by $\psi_j:f(V_j)\to V_j$ the inverse map of $f|_{V_j}$.
Since $f(V_j)$ is not necessarily contained in some $V_i$, we refine the cover $V_j$
by putting $V_{ji}=f(V_i)\cap V_j$ and
set $\psi_{ji}=\psi_j|_{V_{ji}}:V_{ji}\to\psi_j(V_{ji})\subset V_i$.
We choose for each $j$ a $C^r$ partition of the unity $\{\chi_{ji}\}$ on $V_j$
such that each $\chi_{ji}$ is supported in $V_{ji}$.
Finally, we set
$
g_{ji}(a)=\chi_i(\psi_j(a))\cdot\tilde{\chi}_{ji}(a)\cdot g_j(a) 
$
(note that $g_{ji}$ is compactly supported in $V_{ji}$),
and we may rewrite our operator as
$$
\LL_0\phi(a)=\sum_{j,i} g_{ji}(a)\cdot\phi(\psi_{ji}(a)) \, .
\tag{4.3}
$$
The operators $\LL_k$ acting on $k$-forms on $M$ are similarly defined (replacing
the composition
with $\psi_{ji}$ by the pullback). 
We  next choose charts $\pi_j:U_j\to V_j$,
where the $U_j$ are bounded and two-by-two disjoint
open subsets of $\real^n$. We denote
by $\pi$ the
map  $\pi|_{U_j}=\pi_j$ and set
$U_{ji}=\pi_j^{-1}(V_{ji})$.
 
Let now $\varphi$ be a form in $U=\cup_j U_j$. We set
$(\tau\varphi)(a)=\sum_i\chi_i(a)\cdot(\pi_i^{-1})^*\varphi(a)$.
Then, $\tau\pi^*$ is the identity on $C^r$ forms in $M$ (in particular $\pi^*$
is injective).
We now define $\MM_k$
acting on $k$-forms in $U\subset\real^n$ as $\MM_k=\pi^*\LL_k\tau$, that is
$$
\MM_k\varphi(x)=\sum_{ji}g_{ji}(\pi(x))
\chi_i(\psi_{ji}\circ\pi(x))\cdot(\pi_i^{-1}\circ\psi_{ji}\circ\pi)^*\varphi(x) \, .
\tag{4.4}
$$
It is easy to choose $\Omega$, $\psi_\omega$ and $g_\omega$
in order to view
$\MM_k$ as an operator of the form \thetag{2.3}, with $g_\omega$ compactly supported in
$U_\omega$.

\smallskip
{\bf Sobolev spaces}
\smallskip
We shall work  with
the Bessel potential $\JJ_1$,  defined on $L^{q}(\real^n)$ by
(see [St, Chapter ~V.3], $\Gamma(\cdot)$ denotes the Euler
gamma function):
$$
\JJ_1 (\varphi)(x)={1 \over \Gamma(1/2)\sqrt{4\pi}}
\int_{\real^n}
\int_0^\infty e^{-\pi |y|^2/t} e^{-t/4\pi}
t^{(-n+1)/2} {dt \over t} \,
\varphi(x-y) \, dy   \, .
$$

Write $\JJ_\ell$ for the  $\ell$th iterate of $\JJ_1$.
It is a well-known (see [St, V.3.3-3.4])
and important result in the theory
of Sobolev spaces that 
for all $1 < q < n$ and each $\ell \ge 1$, $\JJ_\ell$ is an isomorphism from the Sobolev
space $L^q_m(\real^n)=W^{m,q}(\real^n)$ to the Sobolev space
$L^q_{m+\ell}(\real^n)$
if $1< q < \infty$, $m$ is a nonnegative integer,
and $\ell \ge 1$
[St, V.3.2-3.3]. In particular, $\JJ_1^{-1}
\varphi\in L^q_{m-1}(\real^n)$
if $\varphi \in L^q_m$ and $m \ge 1$. 
We take $\widetilde \JJ_k=\JJ_{r-k}$.

\smallskip

We define the Sobolev spaces $W^{m,p}(M)$ [Ad] of the manifold $M$ using our chart
$\pi:U\to M$ (by compactness, other chart systems yield equivalent norms):
$$
W^{m,p}(M):=\{ \phi:M\to\complex\, |\, \phi\circ\pi \in W^{m,p}(U)\}.
$$
For $\phi\in\AA_{k,W^{m,p}(M)}$, set $||\phi||_{\AA_{k,W^{m,p}(M)}}:= 
||\pi^*u||_{\AA_{k,W^{m,p}(U)}}$.

Clearly, $\LL_k$ is bounded on $\AA_{k,W^{m,p}(M)}$ for all $0<p<1$ and 
$m \le r-1$
(for $k=0$ we can take $r=m$).
By definition, $\pi^*$ is a Banach space isomorphism between $\AA_{k,W^{m,p}(M)}$ and
its image in $\AA_{k,W^{m,p}(U)}$. The spectrum of  $\LL_k$  on $\AA_{k,W^{m,p}(M)}$ 
coincides with its spectrum on the image $\LL_k(\AA_{k,W^{m,p}(M)})$, similarly
for $\MM_k$ and  $\AA_{k,W^{m,p}(U)}$. Also, $\pi^*$ is an isomorphism
between $\LL_k(\AA_{k,W^{m,p}(M)})$ and $\MM_k(\AA_{k,W^{m,p}(U)})$, with inverse
$\tau$, and which conjugates $\MM_k$ and $\LL_k$.
By definition,  $\Det^\flat(\Id-z\LL_k)=\Det^\flat(\Id-z\MM_k)$. We may henceforth concentrate
on the $\MM_k$ acting on $\AA_{k,W^{m,p}(U)}$.

\smallskip
{\bf Checking Axioms 1 and 2--3}
\smallskip

From now on, assume that the data $\psi_\omega$,
$g_\omega$
is $C^{\hat r}$
for some $\hat r\ge nr$ 
and satisfies the hypotheses of Section 2. Additionally, the
$\psi_\omega$ are uniform $\theta$-contractions with
$\theta<1$. We l use the notations of Sections 2 and 3,
in particular the
definitions of $K$ $,K'$ and $\chi_K$, and we assume that $U=K'$.
The Banach spaces
$\BB_{k,t}=W^{rn-k,t}(K')$ clearly satisfy Axiom ~1.
To establish the bounds in Axiom 2  the following estimate will be
instrumental.  We let $\MM_k^+$ be the transfer operator with $g$ replaced
by $|g|$.

\proclaim{Lemma 15 (Lasota-Yorke inequality)} For $0\le k\le n$, $0\le m\le \hat r-1$
($0\le m \le \hat r$ if $k=0$), and
$1 < t < \infty$
there are $C>0$ and a sequence  $\{\rho_\ell\}$
of positive real numbers so that for all
$\ell$  and $\varphi\in W^{m,t}(K')$,
  $$
  ||\MM_k^\ell\varphi||_{m,t}\le \rho_\ell ||\varphi||_{m-1,t}
  + C\cdot\theta^{\ell m }||(\MM_k^+)^{\ell }||_t\cdot ||\varphi||_{m,t}\, .
  \tag{4.5}
  $$
\endproclaim

\demo{Proof of Lemma 15}
The norm $\|\cdot \|_{m,t}$ on $W^{m,t}(K')$ is given by
$
\left( \sum_{|\alpha|\le m} \|D^\alpha \varphi \|_t^t\right)^{1/t}   
$,
where
$\|\cdot \|_t$ is the $L^t(K')$ norm,
$\alpha=(\alpha_1,\dots,\alpha_p)$ is a multi-index of size $|\alpha|=p$
and $D^\alpha=\partial_{\alpha_1}\cdots\partial_{\alpha_p}$.  

We begin with $\ell=1$ and $k=0$. Let $\alpha$ be a multi-index with 
$|\alpha|= m$.
The Leibnitz rule gives
$D^\alpha \MM_0\varphi=\sum_\omega g_\omega D^\alpha(\varphi\circ\psi_\omega) 
+ H_\alpha\varphi$, where $H_\alpha\varphi$ 
is a sum of
terms 
$$
D^{\beta}g_\omega\cdot 
 D^{\beta'} \left ( \varphi  \circ \psi_\omega\right  )\, , 
\tag{4.6}
$$
with  $|\beta|\le |\alpha|$, $|\beta'|\le |\alpha|-1$. The data
being at least $C^m$, all  terms  \thetag{4.6} are bounded by a constant multiple of
$||\varphi||_{m-1,t}$.
Now,  $D^\alpha(\varphi
\circ\psi_\omega)$ is a sum
of, on the one hand, terms involving derivatives of $\varphi$ of order $\le|\alpha|-1 $ 
(which are also bounded by a constant multiple of
$||\varphi||_{m-1,t}$), and, on the other hand, 
terms of the form
($(\psi_\omega)_j$ denotes the $j$-th coordinate of $\psi_\omega$)
$$
 \sum_\omega g_\omega\cdot \bigl(\prod_{i=1}^m\partial_{\alpha_i}
 ((\psi_\omega)_{\gamma_i})\bigr)\cdot
 \left( D^\alpha \varphi\right)\circ\psi_\omega \, .
 \tag{4.7}
$$
The $L^t$ norm to the power 
$t$ of \thetag{4.7}
is $\le\theta^{mt}||\MM_0^+(D^\alpha \varphi)||_t^t$, by the contraction
assumption. Since $\varphi\in W^{m,t}(K')$, 
$||D^\alpha\varphi||_t\le||\varphi||_{m,t}$.
The number of terms  \thetag{4.7}  depends only on $m$ and $n$.
The result for $\ell=1$ follows by summing 
over $\alpha$ with $|\alpha|\le m$
(if $|\alpha|\le m-1$, terms of the form
\thetag{4.7} are also bounded  by a constant multiple of
$||\varphi||_{m-1,t}$).

For $\ell>1$, just use that $\|(\MM_0^\ell)^+\|_t\le \|(\MM_0^+)^\ell\|_t$.

The claims for $k\ge 1$ are obvious by the contraction property of the
pullback.
\qed
\enddemo

\proclaim{Lemma 16 (Axiom 2)} Let $1 < t < \infty$.

\noindent (1)
$\rho_{ess}(\MM_0|_{W^{m,t}(K')})\le \theta^{m} \rho(\MM^+_0|_{L^{t}(K')})$
for all $0\le m\le \hat r$.

\noindent (2) For all $0\le m\le \hat r-1$ we have
$\rho(\MM_k|_{\AA_{k,L^{t}(K')}})\le \theta^{k}\rho(\MM^+_0|_{L^{t}(K')})$ 
and
$$
\rho(\MM_k|_{\AA_{k,W^{m,t}(K')}})
\le \theta^k \rho(\MM^+_0|_{L^{t}(K')})  \, , \,\,
\rho_{ess}(\MM_k|_{\AA_{k,W^{m,t}(K')}})\le \theta^{m+k} \rho(\MM^+_0|_{L^{t}(K')})\, .
$$

\noindent (3)
$\rho(\MM^+_0|_{L^{t}(K')})\le 1 = e^P$ if $g_\omega = |\det D\psi_\omega|$. Otherwise
$\rho(\MM^+_0|_{L^{t}(K')})\le \theta ^{-n/t} e^P$.

The spaces 
$\BB_{k,t}=W^{\hat r-k,t}(K')$ for $k=0, \ldots, n$, together with
the  operators $\MM_k$,
satisfy Axiom ~2 for 
$\widetilde R=\theta^{\hat r}$ 
if  $g_\omega=|\det D\psi_\omega|$, and
(after renormalising so that $e^{P(\log |g|}=1$)  and $\widetilde
R> \theta^{\hat r-(n-1)}$ 
otherwise (if $t$ is not much smaller than $n/(n-1)$).
\endproclaim

\demo{Proof of Lemma 16}
For the first claim,
since $W^{m-1,t}(K')$ is compactly embedded in $W^{m,t}(K')$ [Ad],
we may combine
Lemma~ 15 and the Hennion formula ([He]) to deduce
$
\rho_{ess}(\MM_0|_{\AA_{0,W^{m,t}(K')}})\le \theta^m\rho(\MM^+_0|_{\AA_{0,L^t(K')}})
$.

The  inequalities in the second claim are then obvious, by the action of the pullback.

To prove the third claim, it is more convenient to estimate
the spectral radius of $\LL_0^+$ acting on $L^t(M)$
(see [MS] for an analogous bound in a more specific
situation). The H\"older inequality
for finite sums gives, for $1 < t' < t$ so that $t^{-1}+(t')^{-1}=1$,
and writing $\LL_{0,t'}$ for the operator associated to
$|\det (D\psi_\omega )|\cdot (|g_\omega|/|\det D\psi_\omega|)^{t'}$
$$
\align
&\|(\LL^+_0)^\ell \varphi \|^t_{L^t(M)} = 
 \int_{M} \bigl ( \sum_{\vec \omega \in\Omega^\ell} 
\prod_{j=1}^\ell |g_{\omega_j} \circ\psi^{j-1}_{\vec \omega}| \cdot
|\varphi\circ \psi^\ell_{\vec \omega}| \bigr )^t \, dx\cr
&\quad \le
\int_{M} 
\biggl (  
\sum_{\vec \omega \in\Omega^\ell}  
|\varphi\circ \psi^\ell_{\vec \omega}|^t \cdot|\det D\psi^\ell_{\vec \omega}|
\biggr )
\biggl ( 
\sum_{\vec \omega \in\Omega^\ell} 
\bigl (
\prod_{j=1}^\ell |g_{\omega_j} \circ\psi^{j-1}_{\vec \omega}| \cdot
|\det D\psi^\ell_{\vec \omega}|^{-1/t} \bigr )^{t'}
\biggr )^{t/t'} \, dx
\cr
\allowdisplaybreak
&\quad \le
\bigl ( \int_{M}  |\varphi|^t \, dx \bigr ) \cdot 
\sup_M \biggl ( 
\sum_{\vec \omega \in\Omega^\ell} 
|\det D\psi^\ell_{\vec \omega}| \cdot \bigl (
{\prod_{j=1}^\ell |g_{\omega_j} \circ\psi^{j-1}_{\vec \omega}| \over
|\det D\psi^\ell_{\vec \omega}|} \bigr )^{t'}
\biggr )^{t/t'} \, dx
\cr
&\quad \le
 \|\varphi\|^t_{L^t}\cdot 
\bigl ( \sup_M ((\LL_{0,t'})^\ell  (1)) \bigr )^{t/t'}\, , 
\qquad\forall \ell \, , \forall \varphi \in L^t(M) \, .
\tag{4.8}
\endalign
$$
Taking the $t$-th root of the above inequality
implies 
$$
\rho(\MM_0|_{L^t(K')}) =\rho(\LL_0|_{L^t(M)}) \le 1=e^{P(\log g)} \hbox{ if }
g_\omega=|\det D\psi_\omega| \, .
$$
If 
$g_\omega\ne |\det D\psi_\omega|$, we get 
$$
\rho(\LL^+_0|_{L^t(M)})
\le \lim_{\ell \to \infty}
\bigl (\sup_M \bigl ( 
\sum_{\vec \omega \in\Omega^\ell} 
|\det D\psi^\ell_{\vec \omega}|^{-1+1/t'} \cdot 
\prod_{j=1}^\ell |g_{\omega_j} \circ\psi^{j-1}_{\vec \omega}|  \bigr ) \bigr )^{1/\ell}
\le \theta^{-n/t}\cdot
e^{P(\log |g|)}\, . 
$$

This immediately implies Axiom~2 for the $\MM_k$ and the claimed value of $\widetilde R$.
\qed
\enddemo

The properties in Axiom 3 hold by uniform contraction and bounded distortion.

\smallskip

{\bf Completing the new proof of Ruelle's theorem}
\smallskip

The proofs of Lemmas~15--16 adapt to the
$C^m$ setting, giving for all $f$, $g$ the well-known bound
for all $0\le m \le \hat r-1$ 
($0\le m\le \hat r$ if $k=0$)
$$
\rho_{ess}(\LL_k|_{\AA_{k,C^{m}(M)}})=
\rho_{ess}(\MM_k|_{\AA_{k,C^{m}(\overline{K'})}})\le \theta^{m+k}e^P
 \, .
$$ 
Indeed,
when proving the Lasota-Yorke inequality, we may work with $\MM_0$ acting on $C^0$ instead
of $L^t$ and the extraneous $\theta^{-n/t}$  factor does not appear.
(This also proves the estimates as $t\to \infty$ in Axiom~2(2).)

Clearly, a generalised eigenfunction in $\AA_{k,C^{\hat r-k}(\overline{K'})}$
is in $\AA_{k,W^{\hat r-k,t}(K')}$. The converse is true, working with
eigenfunctions of the dual.
Hence, if $\hat r \ge rn$
and   $f$ and $g$ are $C^{\hat r}$, the result of Theorem ~13 on the Sobolev spaces in fact
implies the statement given above of Ruelle's theorem if $g(y)=1/|\det D f(y)|$
and a weaker statement, replacing $\theta^{-r}e^{-P}$ by 
$\theta^{-(r-1+1/n)}e^{-P}$, for general $C^{\hat r}$ weights.

\enddocument